\documentclass[12pt]{article}
\usepackage{amsfonts,amsmath,latexsym,amssymb,mathrsfs}
\usepackage{graphicx}

\newenvironment{proof}{\noindent{{\bf Proof}:\ }}{\vspace{2ex}}

\newenvironment{remark}{{\nin\bf Remark.}}{{\vspace{2ex}}}
\def\numberlikeadb{\global\def\theequation{\thesection.\arabic{equation}}}
\numberlikeadb
\newtheorem{theorem}{Theorem}[section]
\newtheorem{lemma}[theorem]{Lemma}
\newtheorem{corollary}[theorem]{Corollary}

\newcommand{\RR}{{\mathbb R}}

\def\nat{{\mathbb N}}
\newcommand{\var}{{\mbox{Var}}}

\newcommand{\beas}{\begin{eqnarray*}}
\newcommand{\enas}{\end{eqnarray*}}
\newcommand{\eqs}{\begin{eqnarray*}}
\newcommand{\ens}{\end{eqnarray*}}

\newcommand{\bea}{\begin{eqnarray}}
\newcommand{\eqa}{\begin{eqnarray}}
\newcommand{\ena}{\end{eqnarray}}
\newcommand{\eq}{\begin{equation}}
\newcommand{\en}{\end{equation}}

\newcommand{\proofbox}{\hspace*{\fill}\mbox{$\halmos$}}
\newcommand{\halmos}{\rule{1ex}{1.4ex}}

\def\half{{\textstyle{\frac12}}}

\def\third{{\textstyle{\frac13}}}

\def\eighth{{\textstyle{\frac18}}}
\def\Ref#1{(\ref{#1})}
\def\a{\alpha}
\def\b{\beta}
\def\s{\sigma}
\def\f{\phi}

\def\l{\lambda}

\def\p{\pi}

\def\dtv{d_{TV}}

\def\ep{\hfill $\proofbox$ \bigskip}

\def\re{\RR}
\def\giv{\,|\,}

\def\Po{{\rm Po\,}}
\def\non{\nonumber}
\def\th{\theta}
\def\e{\varepsilon}
\def\m{\mu}
\def\var{{\rm Var\,}}
\def\r{\rho}

\def\t{\tau}
\def\g{\gamma}
\def\h{\eta}
\def\f{\phi}
\def\ps{\psi}

\def\lti{{\lim_{t\to\infty}}}

\def\Bl{\left(}
\def\Br{\right)}
\def\Blm{\Bigl|}
\def\Brm{\Bigr|}
\def\Blb{\left\{}
\def\Brb{\right\}}

\def\Bi{{\rm Bi\,}}

\def\nin{\noindent}

\def\ex{{\mathbb E}}
\def\pr{{\mathbb P}}
\def\msk{\medskip}

\def\ff{{\cal F}}
\def\d{\delta}
\def\z{\zeta}
\def\Eq{\ =\ }
\def\Le{\ \le\ }

\def\ui{^{(1)}}
\def\bs{{\hat\sigma}}

\def\bB{{\widehat B}}
\def\bL{{\widehat L}}

\def\bZ{{\widehat Z}}
\def\bW{{\widehat W}}
\def\bF{{\widehat F}}
\def\sjjn{\sum_{j \in J_N}}
\def\sli{\sum_{l\ge1}}
\def\sjnnl{\sum_{j\in J_{Nl}}}
\def\skl{\sum_{k=1}^l}
\def\fftn{\ff_{\t_N}}
\def\ffptn{\ff^+_{\t_N}}
\def\ss{{\cal S}}
\def\ee{{\cal E}}

\def\BHJ{{Barbour, Holst \& Janson}}
\def\nfe{\lfloor N^{5/8}\rfloor}
\def\nit{\lfloor \sqrt N \rfloor}
\def\nitmi{\lfloor \sqrt N (\m+1)\rfloor}
\def\nith{\lfloor N^{1/3} \rfloor}

\def\iid{independent and identically distributed}
\def\tA{{\widetilde A}}
\def\tL{{\widetilde L}}
\def\tI{{\tilde L}}
\def\hL{{\widehat L}}
\def\hU{{\widehat U}}
\def\MN{{\rm MN}\,}
\def\bt{{\bar\tau}}
\def\ts{{\tilde s}}
\def\hs{{\hat s}}
\def\ul{^{(l)}}
\def\uii{^{(i)}}
\def\uk{^{(k)}}

\def\ignore#1{}

\def\n{\nu}

\def\Giv{\,\Big|\,}

\def\integ{{\mathbb Z}}

\def\Def{\ :=\ }
\def\uu{^{(u)}}
\def\ut{^{(2)}}
\def\uid{^{(i')}}
\def\ttt{{\tilde \t}}
\def\baf{{\hat\ps}}
\def\skid{\sum_{k=1}^d}
\def\slid{\sum_{l=1}^d}
\def\bxi{{\hat\xi}}
\def\skK{\sum_{k=1}^K}

\def\slK{\sum_{l=1}^K}
\def\hz{{\hat \zeta}}
\def\tW{W_*}
\def\tbW{\bW_*}
\def\tbaf{\baf_*}
\def\hmu{{\hat\mu}}
\def\bm{\hmu}
\def\heta{{\hat\h}}
\def\uit{^{(1)^T}}
\def\hZ{{\widehat Z}}
\def\hH{{\widehat H}}
\def\diag{{\mbox{diag}}}
\def\tf{{\tilde f}}
\def\tq{{\tilde q}}

\begin{document}

\title{Approximating the epidemic curve}
\author{
A. D. Barbour\footnote{Institut f\"ur Mathematik, Universit\"at Z\"urich,
Winterthurertrasse 190, CH-8057 Z\"URICH.
Work supported in part by Australian Research Council Grants Nos DP120102728 and DP120102398
\msk}
\ and
G. Reinert\footnote{Department of Statistics,
University of Oxford, 1 South Parks Road, OXFORD OX1 3TG, UK.
GDR was supported in part by EPSRC and BBSRC through OCISB.
}\\
Universit\"at Z\"urich and University of Oxford }

\date{}
\maketitle

\begin{abstract}
Many models of epidemic spread have a common qualitative structure.  The numbers
of infected individuals during the initial stages of an epidemic can be well 
approximated by a branching process, after which the proportion of individuals 
that are susceptible follows a more or less deterministic course.  In this
paper, we show that both of these features are consequences of assuming a
locally branching structure in the models, and that the deterministic course
can itself be determined from the distribution of the limiting random variable
associated with the backward, susceptibility branching process.  Examples considered include
a stochastic version of the Kermack \& McKendrick model, the Reed--Frost model, 
and the Volz configuration model.
\end{abstract}

\nin{{\bf Keywords.}  Epidemics, Reed--Frost, configuration model, deterministic approximation,
branching processes} 

\msk\nin
{{\bf Maths Reviews Classification.}} 92H30, 60K35, 60J85


\section{Introduction}\label{intro}
 \setcounter{equation}{0}
Kermack \& McKendrick's (1927) 
model of the course of an epidemic in a closed population
has proved to be both effective in practice (see for example Brauer~(2005),  
Brauer \& Castillo--Chav\'ez~(2012) p.350, Gupta {\it et al.\/} (2011))
and influential in the theoretical development of epidemic modelling.  Writing
$s(t)$ to denote the {\it density\/} of susceptible individuals in the population
at time~$t$ and $\b(v)$ the infectivity of an individual at time~$v$
after becoming infected, and normalizing the initial population density to be
$s(-\infty)=1$,  the development of~$s$ is given by the equation
\eq\label{K-McK}
  (- Ds(t)) \Eq s(t)\int_0^\infty \b(v)(- Ds(t-v))\,dv.
\en
Here, $Ds$ denotes the derivative of~$s$ with respect to time, and is {\it negative\/}.
The quantity $(-Ds(t))$ is the rate at which the density of susceptibles is being reduced
at time~$t$, and this is just the (density standardized) rate at which infections are 
being made, explaining the integral on the right hand side of~\Ref{K-McK} as the
force of infection at time~$t$.  Dividing both sides of~\Ref{K-McK} by~$s(t)$ and
integrating gives 
\eq\label{deterministic}
     -\log s(t) \Eq \int_0^\infty \b(v)\{1 - s(t-v)\}\,dv  \,.
\en
Note that, if~$s$ satisfies~\Ref{deterministic}, so does any translate~$\ts_h$ defined 
by $\ts_h(t) = s(t+h)$, for any~$h \in \re$.  However, it is shown in Diekmann~(1977) that
there is exactly one solution~$s$ to~\Ref{deterministic} that is non-increasing and non-negative, 
if, for instance, the value of~$s(0) \in (0,1)$ is specified.
Letting $t \to \infty$, \Ref{deterministic} gives the final size equation
\eq\label{final-size}
    - \log s(\infty) \Eq R_0(1 - s(\infty))\,,
\en
where the basic reproduction number~$R_0 := \int_0^\infty \b(v)\,dv$ is the
expected total number of infections made by an infected individual in a susceptible 
population of unit density;  a proportion $1-s(\infty)$ of the population has 
been infected by the end of the epidemic.  Kermack \& McKendrick~(1927) then
deduced their famous {\it threshold theorem\/}, that $s(\infty) < 1$ is
only possible if $R_0 > 1$.

The final size equation can be interpreted more directly, without integrating~\Ref{K-McK},
but at the level of an individual. Rewrite~\Ref{final-size} in the form
\eq\label{expl-final-size}
   s(\infty) \Eq e^{- R_0(1 - s(\infty))},
\en
and recognize $R_0(1-s(\infty))$ as the total integrated force of infection over the
whole course of the epidemic.  The tacit assumption about force of infection at the
level of the individual is that it represents the `instantaneous rate' of infection
of an individual, interpreted in a Markovian sense, so that the probability of an
individual avoiding infection after exposure to an integrated force of infection~$f$ 
should be given by~$e^{-f}$.  Thus the right hand side of~\Ref{expl-final-size}
is the probability that an individual avoids infection throughout the whole course 
of the epidemic, which is exactly the proportion~$s(\infty)$ that remain uninfected
to the end.  

The equation~\Ref{expl-final-size}, with~$s(\infty)$ replaced by the symbol~$q$,
also has a classic interpretation in a branching process context.  It represents
the equation for the extinction probability~$q$ of a branching process starting 
with a single individual, when the number of offspring has the Poisson 
distribution~$\Po(R_0)$ with mean~$R_0$.  
At an individual level, this suggests
a stochastic analogue of the Kermack--McKendrick model, in which an infected
individual makes potentially infectious contacts according to a Poisson process of 
rate~$\b(v)$, where~$v$ represents the time since infection. 
Each such event leads to a new infection, if the individual contacted is
susceptible.  In the early stages of an epidemic, almost all individuals are
still susceptible, and so the early development of the epidemic is well
approximated by a branching process, in which an individual at age~$v$
has (Markovian) birth rate~$\b(v)$.
Branching processes have long been used to approximate the early stages of 
epidemic processes in this way.  The earliest papers are those 
of Kendall~(1956) and Whittle~(1955), and a systematic treatment is given 
in Ball \& Donnelly~(1995).  In particular, the Kermack--McKendrick threshold theorem is replaced by a
stochastic threshold theorem, in which the probability that a large
epidemic takes place, when started by a single infected individual~$K_0$ in an initially susceptible 
population of large size~$N$, is (approximately) $1-q$, thus being positive exactly
when the mean number of offspring, here~$R_0$, exceeds~$1$.

In contrast, for the 
analysis of the final size $N - S(\infty)$, 
where $S(t)$ denotes the number of susceptibles at time $t$, 
the appropriate branching approximation is not at the
beginning of the epidemic, but approximates the process of the contacts potentially
leading to the infection of a randomly chosen individual~$K$ --- see, for
example, Diekmann \& Heesterbeek~(2000), pp.~171--172.  If this backward process of contacts
contains few individuals, as when its branching approximation dies out, then~$K$
is unlikely to become infected, whereas, if it contains many individuals, 
as when the branching approximation never dies out,  
$K$ is almost certain to become infected, if the epidemic is a large one. 
Thus the probability that a randomly chosen~$K$ does not become infected is
approximately~$1$ if the epidemic starting from~$K_0$ is a small one, and
approximately the extinction probability~$q_b$ for the `backward' branching process,
if the epidemic is a large one.  
However,
because of the random choice of~$K$, the probability that~$K$ escapes infection is just
$N^{-1}\ex S(\infty)$.  Hence 
$$
  1 - N^{-1}\ex S(\infty) \ \approx\ 1 - \{q + (1-q)q_b\} \Eq (1-q)(1-q_b),
$$
so that, given that the epidemic starting from~$K_0$ is a large one,
the (mean of the) final proportion of infected individuals
is close to~$(1-q_b)$.  As it happens, for the stochastic Kermack--McKendrick model 
described above, the forward and backward branching processes are the same, so that
$q_b=q$, and~\Ref{expl-final-size} is still the relevant equation for determining
the final outcome of the epidemic, with $s(\infty)$ replaced by~$q_b$. 
Thus,  
in the deterministic model, a large epidemic is certain, and
the proportion of the population that is infected is $(1-q_b)$. In the stochastic
model, a large epidemic occurs only with probability approximately $(1-q)$, in which
case a proportion of approximately $(1-q_b)$ of the individuals are infected, and
on the complementary event there is only a tiny outbreak involving a negligible
proportion of infected individuals.  However, if the epidemic were started with~$I > 1$
individuals, the probability of a large outbreak, again leading to 
a proportion of approximately 
$(1-q_b)$ of the individuals being infected, increases to $(1 - q^I)$, and is thus
nearly a certain event if~$I$ is at all large.

In this paper, we use analogous ideas to show that, under appropriate
assumptions, the whole course of the
stochastic epidemic is determined by the analysis of the two branching processes,
forward and backward.  There is an initial phase, approximated as usual by the
forward branching process.  If this branching process does not become extinct, it
settles to an essentially deterministic course of exponential growth, after a random delay 
that results from the initial random development of the branching process.  After the
point at which the forward branching process ceases to be a good approximation, 
the proportion of susceptibles in the epidemic process follows an
almost deterministic development, which can be expressed in terms of properties of the backward
branching process.  
One of the consequences of this is to show that the Markovian stochastic
interpretation of the instantaneous force of infection, which is implicit in the derivation
of the deterministic Kermack--McKendrick equation~\Ref{K-McK}, is not actually necessary
to justify the equation; we prove that~\Ref{K-McK} holds as a faithful approximation in a much
wider class of models.

We illustrate the approach for the Reed--Frost discrete generation epidemic model
in a population of size~$N$.
Let the probability of an infected individual infecting a given susceptible
be $p = \m/N$. Then the approximating Galton-Watson forward branching process has offspring
distribution~$\Po(\m)$ (and $R_0=\m$); we take $\m > 1$.  
After~$n$ time units, the number of  individuals alive in
the branching process is $Z_n \sim W\m^n$ and the total number of individuals that were alive 
in previous generations is approximately $W\m^n/(\m-1)$, where~$W$ is the a.s.\ limit of~$Z_n \m^{-n}$.  Take
\[
    n \Eq n(N) \Def \lfloor \half\log N/\log \m \rfloor,
\]
so that $\m^n = \th_N N^{1/2}$, where $1\le \th_N < \m$, and suppose that~$W > 0$.
Label those that have died in chronological order, with labels drawn independently 
and at random from~$[N] := \{1,2,\ldots,N\}$.  Mark any whose labels have been
used before, and all of their descendants, as `ghosts'.  There are only few
marked, and those that are unmarked are
the individuals that have been infected before time~$n$ in the epidemic. Let the 
set of labels used be denoted by~$L_N$; its size is small compared to~$N$.

Now, starting from a randomly chosen individual, take an independent realization of the 
reversed branching process --- in this model, it has
the same law as the forward process --- and run it for $n(N) + r$ generations,
after which there have been approximately $\bW \m^{n+r+1}/(\m-1)$ individuals born in total,
where~$\bW$ is the corresponding realization of the limit random variable, and
is independent of~$W$.  Label these individuals in chronological order at random 
from~$[N]\setminus L_N$, and again mark the (few) ghosts; let the set of labels be~$L^b_N$, 
and denote by~$K$ the label of the initial individual. Do the same for the 
individuals alive in generation~$n$ of the forward process, and call this set~$L^f_N$. If  
$L^b_N \cap L^f_N \ne \emptyset$, and an element of the intersection is a non-ghost,
we can construct a chain of infection to it from the initial individual in the epidemic,  
and a chain going from it to~$K$, 
giving a chain of infection from the start of the epidemic to~$K$. Conversely, any chain 
of infection from the start of the epidemic to~$K$ must pass through 
a non-ghost element of $L^b_N \cap L^f_N$.
Thus there is no chain of infection from the start of the epidemic to~$K$ exactly
when $L^b_N \cap L^f_N$ is empty or contains only ghosts; the event that
$L^b_N \cap L^f_N$ is non-empty but contains only ghosts has only small probability.

Now, given~$Z_n$ and the realization of the 
backward branching process, the mean number of intersections between $L^b_N$ and~$L^f_N$ is 
close to $N^{-1}Z_n\,\bW \m^{n+r+1}/(\m-1)$, and hence, using a Poisson approximation, 
the probability of the intersection being empty is close to 
\eqs
    \exp\{-N^{-1}Z_n\,\bW \m^{n+r+1}/(\m-1)\} &=& \exp\{-N^{-1/2}Z_n\,\bW \th_N \m^{r+1}/(\m-1)\}.
\ens
It is now easy to convert this result into the statement
\eqs
    \lefteqn{\pr[K \mbox{ has escaped infection until generation } 2n+r \giv \ff_n]} \\
     &&\sim\  \ex\bigl\{ \exp\{-N^{-1/2}Z_n\,\bW \th_N \m^{r+1}/(\m-1)\} \giv \ff_n \bigr\},
\ens
where~$K$ is a randomly chosen label from all of~$[N]$ and~$\ff_n$ denotes $\s(Z_l,\,0\le l\le n)$; 
in other words, still with $n = n(N)$,
\[
    \ex\{N^{-1}S_N(2n+r) \giv \ff_n\} \ \sim\ \ps(N^{-1/2}Z_n  \th_N \m^{r+1}/(\m-1)),
\]
where $S_N(t)$ denotes the number of susceptibles in the epidemic at generation~$t$ and
$\ps(\th) := \ex\{e^{-\th\bW}\}$.  

But now, for two independently randomly chosen individuals
$K$ and~$K'$,
\eqs
    \lefteqn{\ex\{(N^{-1}S_N(2n+r))^2 \giv \ff_n\}}\\ 
     &&\Eq \pr[\mbox{both }K\mbox{ and }K' \mbox{ have escaped infection until generation } 2n+r \giv \ff_n]
\ens
can be approximated in exactly the same way; since there is little overlap between the labels 
assigned to the backward branching processes starting from $K$ and~$K'$, it is easy to deduce that
$$
    \ex\{(N^{-1}S_N(2n+r))^2 \giv \ff_n\} \ \sim\ \{\ps(N^{-1/2}Z_n  \th_N \m^{r+1}/(\m-1))\}^2
$$
also, implying that $\var\{N^{-1}S_N(2n+r) \giv \ff_n\} \sim 0$.  
Writing $W=\lim_{m \rightarrow \infty} Z_m \mu^{-m}$, we note that 
$Z_n = Z_{n(N)} \sim W \mu^{n(N)} = \theta_N N^{\frac12} W$;
this implies that, for any $\e > 0$ and any $r\in\integ$,
\eq\label{Reed-Frost}
    \lim_{N\to\infty} \pr[|N^{-1}S_N(2n(N)+r) - \ps(W \th_N^2 \m^{r+1}/(\m-1))| > \e]
     \Eq 0.
\en
The quantity $\ps(W \th_N^2 \m^{r+1}/(\m-1))$ is random only through the presence of~$W$.
By time~$n(N)$ the quantity $W$ is essentially determined, and is the same for 
all~$r\in\integ$.  If $W = 0$, the above approximation is by $\ps(0) = 1$ for all~$r$,
indicating that only a small epidemic occurs; the assumption $\m > 1$ merely ensures that
$\pr[W>0] > 0$, so that a large epidemic is indeed possible.

If $W>0$, one could describe the approximation slightly differently.
The values of $N^{-1}S_N(2n(N)+r)$ for $r\in\integ$ are then approximated by a 
discrete subset of points on the continuous deterministic curve 
$u \mapsto \ps(\m^{u+1}/(\m-1))$, namely those with~$u$ of the form 
$r + \{\log W + 2\log\th_N\}/\log\m$ for $r\in\integ$.
Thus randomness appears only as a time shift in the lattice
of integer spaced points along the continuous deterministic path that are used for the 
approximation to the discrete time process.  Note also that the times~$l$ at which
$N^{-1}S_N(l)$ is not close either to~$0$ or to~$1$ are within~$O(1)$ of $\log N/\log \m$;
the development of the epidemic is slow until almost time $\log N/\log \m$, and then
runs its course over comparatively few time steps.
  
  In what follows, we shall make these arguments precise, but
for processes with non-lattice offspring distributions in continuous time.  The phenomena
associated with
discretization disappear, giving a neater result, but connecting the forward and backward branching
processes becomes more delicate. Our analogue of~\Ref{Reed-Frost} is proved in Theorem~\Ref{main},
under some fairly mild assumptions on the individual point processes of infection that
include the stochastic Kermack--McKendrick model described above for many choices of the infectivity 
function~$\b$.   It establishes that 
\eq\label{K-McK-limit}
    \lim_{N\to\infty} \pr[\sup_u |N^{-1}S_N(\l^{-1}\{\log N -\log W + u\}) - \hs(u)| > \e] 
     \Eq 0,
\en
for a deterministic function~$\hs$, whenever $W>0$; here, $\l$ is the Malthusian parameter 
(assumed positive) and~$W$ the limiting
random variable for the associated forward branching process, and~$\hs$ is determined by the
properties of the associated backward branching process.
The methods that we use have quite general application, and have already been exploited 
in Barbour \& Reinert~(2012) in the context of the Aldous(2010) gossip process and of the 
Moore \& Newman~(1999) small world model.  

The key ingredients that make the proofs go through
are the branching nature of the forward and backward processes, and their exponential
growth and stability properties.  These are also shared, for instance, by their 
multitype analogues. We give a multitype analogue of~\Ref{K-McK-limit} in Section~\ref{multitype},
and discuss a configuration model in Section~\ref{configuration}.

\section{The single type model}\label{basic}
 \setcounter{equation}{0}
\subsection{The branching processes}\label{branching}
We begin by considering an epidemic in a closed population of~$N$ individuals, 
where~$N$ is to be thought of as large, that evolves
according to the following scheme.  Each individual~$i$, $1\le i\le n$, is
equipped with a potential infection history, in the form of a realization
of a point process~$\xi_i$ on $(0,\infty)$.  If~$i$ becomes infected at 
time~$\s(i) < \infty$, it makes infectious contacts with other individuals
at times $\s(i,j) := \s(i) + \t(i,j)$, where $0 < \t(i,1) \le \t(i,2) \le \cdots$
denote the times of the events of~$\xi_i$ and $\n(i) := \xi_i(\re_+) < \infty$ 
their number; if required, $\xi_i$ can be augmented by a time~$\t^r(i) \ge \t(i,\n(i))$,
indicating that~$i$ is removed from the infectious state at time $\s(i) + \t^r(i)$.
The individuals contacted are chosen independently at random from~$[N]$, and an
infectious contact only results in the individual contacted becoming infected if
they have not previously been contacted.  The epidemic begins with individual~$i_1$
becoming infected at time~$\s(i_1)=0$.  After the \hbox{$r$-th} individual~$i_r$ has become
infected at time~$\s(i_r)$, and if $r<N$, then potential infectious contacts occur
at the times $\s(i_r) + v_r(j)$, $1 \le j \le |V_r|$, where the $v_r(j)$ are the elements
of
\[
    V_r \Def \{\s(i_l,j) - \s(i_r), \,1\le j\le \n(i_l),\,1\le l \le r\} \cap (0,\infty),
\]
arranged in non-decreasing order, and the labels of the individuals to be contacted are given by
$I_r(j)$, $j\ge1$, chosen independently and uniformly on~$[N]$.
Defining the index 
$j_*(r) := \min\bigl\{1\le j\le |V_r|\colon\,I_r(j) \notin \{i_1,i_2,\ldots,i_r\} \bigr\}$,
then
\[
     i_{r+1} := I_r(j_*(r)) \quad \mbox{and}\quad \s(i_{r+1}) \Eq \s(i_r) + v_r(j_*(r)),
\]
unless there is no such index~$j_*(r)$, in which case the epidemic stops.
It is assumed that $(\xi_i,\,1\le i\le N)$ are independent and identically distributed.

If the labelling were ignored, and $j_*(r)$ were taken to be~$1$ for each~$r\ge1$, and if the
\hbox{$r$-th} infected individual were assigned infection history~$\xi_r'$, with the $(\xi_r',\,r\ge1)$
independent and identically distributed, then the resulting process would be a Crump--Mode--Jagers
branching process~$Z$.  Indeed, if the~$\xi_r'$ are distributed in the same way as~$\xi_1$, 
the paths of the branching and  epidemic processes (neglecting the labelling) can be coupled
so as to agree exactly until $\r := \min\{r\ge0;\, j_*(r) \ge 2\}$ (Ball 1983, Ball \& Donnelly, 1995),
with the epidemic process recoverable from the branching process by adding labelling, and by
marking as `ghosts' individuals infected in the branching process but not in the epidemic process --- 
$(j_*(r)-1)$ such infections occur whenever $j_*(r) \ge 2$ --- together with the individuals
in the branching process that are descended from such individuals.  We shall make substantial
use of this coupling, but only up to times where there have typically been relatively few ghosts
created.

  We shall make the following assumptions on  the distribution of~$\xi_1$ of the above 
Crump--Mode--Jagers branching process. Let $p_j := \pr[\n(1) = j]$ and $\mu = \ex \n(1)$; denote 
the relative intensity measure of~$\xi_1$ by 
\eq\label{G-def}
      G(dt) \Def \m^{-1}\ex \xi_1(dt).
\en

\nin{\bf Assumptions} 
\begin{enumerate}
\item We assume that the branching process is supercritical, and that
$$
     1\ <\ \mu\ <\ \infty;\qquad m_2 \Def \ex\n(1)^2\ <\ \infty.
$$
Let  $\l>0$ denote the Malthusian parameter of the branching process, satisfying
\eq\label{malthus}
     \ex\Bl \int_0^\infty e^{-\l t}\,\xi_i(dt) \Br \Eq 1.
\en
The existence of $\l > 0$ follows from Jagers~(1975), Theorem 6.3.3, pp.131--2.
We write
\eq\label{m-star-def}
    m_* \Def \mu \l \int_0^\infty t e^{-\l t}\,G(dt) \ <\ \infty;
\en
then $m_*/\l$ represents the mean age at child bearing (Jagers~(1989), p.195).
\item The intensity measure $G$ is non-lattice and has finite second moment. 
The support of~$G$ is a finite or semi-infinite open interval $(a,b)$, and
$G(A) \ge \int_A g(x)\,dx$ for any $A \subset (a,b)$, for some continuous positive density~$g$.
If $b = \infty$, then also $g(x) \ge kx^{-\g}$ for all $x \ge x_0$, for some $x_0 > a$, $k > 0$
and $\g > 3$.
\end{enumerate} 

\medskip
\nin{\bf Remark.}\
Strictly speaking, the epidemic might be better modelled by assuming that the labels
assigned to the individuals infected by any given individual~$i$ are chosen at random
\em{without replacement} from the labels \em{excluding}~$i$, and indeed that the number
infected by a single individual cannot exceed~$N-1$.  However, under the assumption
that $m_2 < \infty$, the total variation distance between this distribution of labels
and that being assumed here is at most $\half N^{-1}(m_2+\m)$.  Since we need only to consider
the offspring of at most~$N^{5/8}$ individuals in our calculations, any difference
between the results of the two models occurs with probability of order at most $O(N^{-3/8})$,
and does not affect the results proved in this paper.

\medskip

Letting the infection times in
the branching process be denoted by $(\s'(r),\,r\ge1)$, and writing
\eq\label{B-def}
    B'(t) \Def \max\{r\colon\, \s'(r) \le t\}
\en
for the number of births that have occurred in the branching process by time~$t$, it
follows that $W(t) := B'(t)e^{-\l t} \to W$ a.s.\ for a non-negative random variable~$W$,
(Nerman~(1981), Theorem 5.4), and also that
$\{W > 0\} = \{\lti B'(t) =  \infty\}$ a.s. (see (3.10) in Nerman~(1981)). 
From Corollary 5.6 in Nerman~(1981), and the fact that pointwise convergence to a continuous 
limit of non-decreasing bounded functions on $[0,\infty]$ is always uniform (Jagers~(1975), p.170), 
it also follows that the statistics of the set
\[
    V'(t) \Def \{\s'(l) + \t'(l,j) - t,\,1\le j\le \n'(l),\,1\le l\le r\} \cap (0,\infty),
\]
where $\t'(l,j)$ denotes the $j$-th point of~$\xi_l'$ and $\n'(l) := \xi'_l(\re_+)$,
converge in distribution, as $t\to\infty$, in the sense that, on $\{W > 0\}$,
\eq\label{F_e-def1}
    \lim_{t\to\infty}\sup_{s\ge0}\bigl|\,(|V'(t) \cap (0,s]| / |V'(t)|) - F(s)\bigr| \Eq 0\quad  {\rm a.s.}\,.
\en
Here~$F$ is the distribution function on~$\re_+$ given by 
\eq\label{F_e-def2}
    1 - F(s) \Def \frac{\m}{\m-1}\, \int_s^\infty (1 - e^{-\l(u-s)})\,G(du).
\en

For the epidemic, the corresponding quantities depend on the choice of~$N$, because of
the role played by the labelling in its definition. We define
\[
    B_N(t) \Def \max\{r\colon\, \s(i_r) \le t\}
\]
and, in the natural notation, 
\[
    V_N(t) \Def V_{B_N(t)} + \s(i_{B_N(t)}) - t.
\]
Provided that~$t$ is not too large, $B_N(t)$ is not very much smaller than~$B'(t)$, and
$|V'(t) \setminus V_N(t)|$ is also relatively small.   This is the case if we take
\eq\label{t_N-def}
    t \Eq t_N(u) \Def \l^{-1}(\half\log N + u),
\en
for any fixed $u>0$, since then $B'(t_{N}(u)) \sim W e^u \sqrt N$, and hence the number of
indices of~$[N]$ chosen more than once in the construction of the epidemic up to this time has mean 
\[
    N^{-1}\binom{B'(t_{N}(u))}2  \ \sim\ \half W^2 e^{2u},
\]
of relative order~$O(N^{-1/2})$ when compared to~$B'(t_{N}(u))$ as~$N$ becomes large; this
observation is made precise later.

We now suppose that $W > 0$, and that the branching and epidemic 
processes have been coupled as described above up to the time~$\t_N := \t(B',\nit)$, 
where $\t(B',r) := \inf\{t>0\colon\,B'(t) \ge r\}$ for any $r>0$.
We denote by~$\fftn$ the corresponding $\s$-field, including the information in  
the sets $V'(\t_N)$ and~$V_N(\t_N)$, but not that of the
labels that are to be assigned to them for the epidemic process.
Since $B'(t)e^{-\l t} \to W$ a.s.\ as $t\to\infty$, it follows that
$B'(t-)/B'(t) \to 1$ a.s.\ also, and hence that $\lim_{N\to\infty}N^{-1/2}B'(\t_N) = 1$
a.s.\ as $N\to\infty$.  Thus
\[
   \t_N \Eq \l^{-1}\{\log B'(\t_N) - \log W(\t_N)\} \ 
     \sim\ \l^{-1}\{\half\log N - \log W\} 
\]
as $N\to\infty$. Note that $B'(\t_N) = \nit$ if $G$ is absolutely continuous with respect to 
Lebesgue measure. 

We now examine whether, and
if so when, a randomly chosen individual $K \in [N]$ becomes infected.  To do so, we begin by
writing 
\eq\label{JN-def}
    J_N \Def [N] \setminus \{i_r,\,1\le r\le B'(\t_N)\} 
\en
to denote the set of indices that have not been used in the definition of the epidemic up to
time~$\t_N$, and we set
\eq\label{edf-defs}
  J_{Nl} \Def \{j \in J_N\colon\, \nu(j)= l\},\quad M_{Nl} \Def |J_{Nl}|
  \quad \mbox{and} \quad M_N \Def \sjjn \n(j) \Eq \sli lM_{Nl}.
\en
We then let
\eq\label{edf-defs-2}
  G_{Nl,k}(x) \Def \frac1{M_{Nl}} \sjnnl I[\t(j,k) \le x] 
\en
denote the empirical distribution function of the times of the $k$-th in order potential infections 
of individuals that have~$l$ such in total, and write
\eq\label{edf-defs-3}
   G_N(x) \Def \frac1{M_{N}} \sli M_{Nl} \skl G_{Nl,k}(x) \Eq  \frac1{M_{N}} \sjjn \xi_j(0,x]
\en
for the overall empirical distribution of the infection times of individuals in~$J_N$.
We introduce the $\s$-field 
\eq \label{ftn-def}
\ffptn = \fftn \bigvee 
\s(\{\t(j,k),\,1\le k\le \n(j), j\in J_N\}).
\en 

If $K \in [N] \setminus J_N$, it has already been infected during the epidemic process before
time~$\t_N$; the conditional probability of this occurring is $\z_N :=  N^{-1}B'(\t_N)$, 
and this is small.  If not, 
it can only have been infected if there is a chain of infection
running backwards from~$K$ to one of the $|V_N(\t_N)|$ individuals in~$J_N$ that were 
infected by individuals in~$[N] \setminus J_N$, but at times after~$\t_N$.  Now
the~$M_N$ infection events originating from individuals in~$J_N$ are directed at
independently and randomly chosen individuals in~$[N]$.  Hence, $K$ is potentially
directly infected as a result of a set of $\Bi(M_N,1/N)$--many events; 
the individuals that infect~$K$ (its generation~$1$ predecessors) were
themselves infected at times preceding the infection of~$K$ by amounts realized through a
Bernoulli$(1/N)$ thinning of the set of $M_N$ times $\{\t(j,k),\,1\le k\le \n(j),\,j\in J_N\}$.
This procedure can be iterated to determine the predecessors in successive generations,
with duplicate choices of a pair $(j,k)$ leading to `ghosts', as before.  In this way,
the susceptibility process, consisting of the chains of potential infection leading to~$K$, 
can be generated from a branching
process~$\bZ_N$ with numbers of offspring having a binomial $\Bi(M_N,1/N)$ distribution,
and occurring at times sampled independently from~$G_N$. 

For the purposes of asymptotics, it is inconvenient to have this branching process
dependent on~$N$.  With some associated error, it can be replaced with a branching
process~$\bZ$ that has a Poisson $\Po(\m)$ offspring distribution, noting that 
\eq\label{mu-def}
    \m \Def \sli l p_l \ \approx\ N^{-1}M_N \ \approx\ \frac{M_N}{|J_N|},
\en
with the birth times independently sampled from the distribution~$G$ defined in~\Ref{G-def}. 
Note that we can write 
\eq\label{G-def-2}
   G = \frac1\m \sli p_l \sum_{k=1}^l G_{lk},
\en
where $G_{lk}$ is the distribution function of the time of the $k$-th event in~$\xi_1$,
conditional on~$\n(1)=l$.  For this branching process, we can define
$\bB(t)$ to be the number of births up to time~$t$, and conclude that, under our assumptions, by 
Theorem 5.4 and (3.10) of Nerman~(1981),
\eq\label{what-cvgce} 
   \bB(t) e^{-\l t}\ \to\ \bW \ \mbox{a.s.}\,, 
\en
for a random variable~$\bW$ that satisfies
$\{\bW > 0\} = \{\lti\bB(t) = \infty\}$ a.s.  Furthermore, letting
\[
   A(t) \Def \{a_t(r)\colon\,1\le r \le \bB(t)\},
\]
where $a_t(r) := t-\bs(r)$ is the age at time~$t$  of the $r$-th individual, it
also follows that, on~$\{\bW > 0\}$, by Corollary 5.6 in Nerman~(1981) together with the 
observation from p.170 of Jagers~(1975),
\eq\label{expl-dist}
  \lim_{t \rightarrow \infty}  \sup_{s\ge0}\,
   \Blm \frac1{\bB(t)} \sum_{r=1}^{\bB(t)} I[a_t(r) \le s] - (1 - e^{-\l s}) \Brm \ =\ 0\quad {\rm a.s.}\,.
\en
Note that, for any $\f \ge 0$, 
\[
   \int_0^\infty e^{-\f t} \m\, G(dt) \Eq \ex\Bl \int_0^\infty e^{-\f t}\,\xi_i(dt) \Br,
\]
so that the branching processes $Z$ and~$\bZ$ indeed have the same Malthusian parameter~$\l$.
We consider this branching process run until time~$t_N(u)$ as in~\Ref{t_N-def}, 
and we show in the next section that it represents a good enough approximation to the process 
of chains of potential infection to~$K$. 

 Finally, we assign labels from~$J_N$
independently and at random to the individuals in the set~$U_N$, whose birth 
times are the elements of~$\t_N + V_N(\t_N)$ --- these are the birth times
in the forward epidemic process that have been determined by time~$\t_N$, but have not
occurred by then --- and also to the set~$\hU_N(u)$ composed of the distinct
individuals among the $\bB(t_{N}(u))$ that are born before~$t_{N}(u)$ in the reverse process.
If the same label is chosen for an individual in~$U_N$, having birth
time $\t_N + v_l$, for some $v_l \in V_N(\t_N)$, and for an individual 
in~$\hU_N(u)$, with birth time $\bs(r) \le t_{N}(u)$, then there is a chain 
of infection to~$K$ of length close to
\eqs
   \t_N + v_l + \bs(r) 
   &=& \l^{-1}\{\log \nit + \half\log N - \log W(\t_N) + u \}  + v_l -a_{t_{N}(u)}(r)\\
   &\sim& \l^{-1}\{\log N - \log W  + u\} + v_l -a_{t_{N}(u)}(r);
\ens
the actual length is $\t_N + v_l + \bs_N(r)$, where $\bs_N(r)$ is the birth time in the~$\bZ_N$
process. 
If, for any such pair, $v_l \le  a_{t_{N}(u)}(r)$, so that the length of the chain
of infection is no greater than $\l^{-1}(\log \nit + \half\log N - \log W(\t_N) + u)$, 
and if the $r$-individual is not a ghost,
then~$K$ is infected before this time; that is, approximately, before time 
$\l^{-1}(\log N - \log W + u)$.

\subsection{Approximating~$\bZ_N$ by~$\bZ$}\label{one-type}
The first step to be justified is that the branching process~$\bZ_N$ with offspring numbers
distributed according to the binomial $\Bi(M_N,1/N)$ distribution and with ages at birth 
independently sampled from~$G_N$, as in \Ref{edf-defs} and~\Ref{edf-defs-3}, 
can be replaced in our considerations by the process~$\bZ$, in which the offspring numbers 
have the Poisson $\Po(\m)$ distribution 
and the ages are sampled independently from~$G$, as in \Ref{mu-def} and~\Ref{G-def}.
We begin by showing that the two constructions lead to the same offspring numbers, with
high probability conditional on~$\ffptn \cap A_N$, at least until the first~$\nfe$ sets 
of progeny have been 
sampled; here, $A_N \in \ffptn$ is a suitably chosen event, 
whose complement has small probability.

\begin{lemma}\label{backward-progeny}
Let 
\eq\label{A_N-def}
  A_N \Def \{|N^{-1}M_N - (1-\z_N)\m| \le N^{-7/16}\} \cap \{\z_N \le N^{-1/2}(\mu + 1)\};
\en 
then $\pr[A_N^c] = O(N^{-1/8})$.
On~$A_N$, it is possible to construct realizations of $\bZ_N$ and~$\bZ$ on the same probability
space, in such a way that the numbers of offspring in the first~$\nfe$ sets of progeny
in the two processes are identical with conditional probability $1 - O(N^{-1/8})$.
\end{lemma}

\proof
We begin by noting from~\Ref{edf-defs} that $M_N :=  \sjjn \n(j)$ is a sum of $N - B'(\t_N)$ 
\iid\ random variables with mean~$\m$ and finite variance.  Hence, by Chebyshev's
inequality,
\eq\label{mu-ratio-bnd}
   \pr[|N^{-1}M_N - (1-\z_N)\m| > N^{-7/16}] \Le N^{-1}\ex\{\n(1)^2\}\,N^{7/8} \Eq O(N^{-1/8}).
\en
Then observe that
\eq\label{B'(t_N)-bnd}
   B'(\t_N) \Le B'(0) + \sum_{j=1}^{\nit-1} X_j,
\en
where $X_j$ denotes the number of offspring of the $j$-th born individual (randomly
ordered in the case of simultaneous births).  Hence, with $B'(0)=1$ and since $\m>1$,
\[
   \pr[\z_N \ge N^{-1/2}(\m+1)] \Le \pr[B'(\t_N) \ge 1 + (\nit-1)\m + \sqrt N] \Le N^{-1/2}\var\n(1),
\]
by Chebyshev's inequality.

Now the total variation distance between $\Bi(M_N,1/N)$ and~$\Po(M_N/N)$ is at
most $1/N$ (\BHJ~(1992), (1.23)), so that branching processes with these two offspring
distributions can be coupled so as to agree until after~$\nfe$ sets of progeny have been 
sampled with failure probability of at most~$N^{-3/8}$.   Then, by considering the
likelihood ratio, $r$ independent samples 
from Poisson distributions with means $\m$ and~$\m'$ can be distinguished with probability
at most $\dtv(\Po(r\m),\Po(r\m')) \le r|\m-\m'|/\sqrt{r\m}$; , see for example Barbour, 
Holst \& Janson~(1992), Theorem I.1.C.  
Hence, if $|N^{-1}M_N- \m| \le N^{-7/16} + \m\z_N$ and $\z_N < N^{-1/2}(\m+1)$,
$\nfe$ samples from $\Po(\m)$ and from~$\Po(M_N/N)$ can be coupled so as to be identical, 
except on an event of probability of order~$O(N^{-1/8})$.  This proves the lemma.  
\ep

We now proceed to the comparison between the age distributions $G_N$ and~$G$.  We assume
henceforth that $N \ge n_1$, where
\eq\label{n1-def}
   n_1 \Def \lceil 4(1+ \mu)^2 \rceil,
\en
so that, on~$A_N$, $\z_N \le \half$, and thus $M_N \ge \half N\m$ if $N\ge n_1$.
Recall the $\s$-field $\ffptn$ from \eqref{ftn-def}.

\begin{lemma}\label{G_N-G}
If $N\ge n_1$, there is an event
$A^*_N \in \ffptn$ having $\pr[(A^*_N)^c] = O(N^{-1/8})$ such that,
for suitably chosen $\e_N = O(N^{-1/6})$, we have
\[
   \pr[|G_N^{-1}(U) - G^{-1}(U)| > \ps_N \giv A_N^*] \Le \h_N,  
\]
where $\h_N := \ps_N + 2\e_N$ and $\ps_N^2 := 2\e_N G^{-1}(1-\e_N)$, and 
where~$U \sim {\rm U}[0,1]$. Note that $\ps_N + \h_N = O(N^{-1/24})$ 
as $N\to\infty$ if~$G$ has finite second moment.  
\end{lemma}

\proof
We begin by using the Dvoretzky--Kiefer--Wolfowitz inequality, in the form given by 
Massart~(1990), which shows that
\[
   \pr[\sqrt{M_{Nl}}\sup_x |G_{Nl,k}(x) - G_{lk}(x)| > z] \Le 2e^{-2z^2}
\]
for any $z > \sqrt{\half\log2}$ and any $k,l$.  Taking $z_{N} := \sqrt{2\log N}$, it follows
that
\eq\label{DKW}
   \pr[(A_{Nl,k}^1)^c] \Le 2N^{-4}
\en
for each~$l,k$, where 
$A^1_{Nl,k} := \{\sqrt{M_{Nl}}\sup_x |G_{Nl,k}(x) - G_{lk}(x)| \le z_{N}\} \in \ffptn$.  
Observe that, for all~$x$,
\eqa
   |G_N(x) - G(x)| &\le& \sum_{l=1}^{\nith} \skl \Blb \frac{M_{Nl}}{M_N} |G_{Nl,k}(x) - G_{lk}(x)| 
           + |M_N^{-1}M_{Nl} - p_l| G_{lk}(x) \Brb  \non \\
      &&\mbox{}\qquad  + \frac1{M_N} \sum_{l > \nith} lM_{Nl} + \sum_{l > \nith} lp_l. \label{G-diff}
\ena
Now, by the Chernoff inequalities (Theorem~2.3 in McDiarmid~(1998)), we have
\eq\label{MNl-pl-bnd}
    \pr[(A^2_{Nl})^c] \Le N^{-3},\qquad l\ge0,
\en
where $A^2_{Nl} := \bigl\{\bigl|M_{Nl} - |J_N| p_l \bigr| \le 4\log N(1 \vee \sqrt{Np_l})\bigr\} \in \ffptn$,
and, by Markov's inequality,
\eq\label{lMNl-sum-bnd}
    \pr\bigl[(A^3_N)^c \bigr] \Le N^{1/6} \sum_{l > \nith} lp_l
            \Le N^{-1/6}\ex\{\n(1)^2\},
\en
where $A^3_N := \{\sum_{l > \nith} lM_{Nl} \le N^{5/6}\} \in \ffptn$.

Define
\[
  A^*_N \Def  A_N \cap \Blb \bigcap_{l=1}^{\nith}\bigcap_{k=1}^l A^1_{Nl,k} \Brb
       \cap \Blb \bigcap_{l=1}^{\nith} A^2_{Nl} \Brb \cap A^3_N \,;
\]
then $\pr[(A^*_N)^c] = O(N^{-1/8})$,
by Lemma~\ref{backward-progeny}, \Ref{DKW}, \Ref{MNl-pl-bnd} and~\Ref{lMNl-sum-bnd}.
On~$A^*_N$, from~\Ref{G-diff}, for all~$x\ge0$, we have
\eqs
    |G_N(x) - G(x)| &\le& \sum_{l=1}^{\nith} l\Blb \frac{\sqrt{M_{Nl}}}{M_N}\sqrt{2\log N}
         + 4M_N^{-1}\log N\{1 \vee \sqrt{Np_l}\} \Brb \\
         && + \left| 1 - \mu \frac{|J_N| }{M_N} \right| +\frac{N^{5/6}}{M_N} +N^{-1/3}\ex\{\n(1)^2\}\\
    &=:& \e_N \Eq O(N^{-1/6}).
\ens 
To justify the order of the bound, note first that, from \Ref{MNl-pl-bnd},  on~$A^*_N$, 
\[
    M_{Nl} \Le \left\{\begin{array}{ll}
                       8\log N, &Np_l < 1;\\[1ex]
                       8\log N\sqrt{Np_l}, &1 \le Np_l < \{4\log N\}^2;\\[1ex]
                       2Np_l, &Np_l > \{4\log N\}^2,
                      \end{array}  \right.
\]
and then $\sum_{l \le \nith}l \le N^{2/3}$, $M_N \ge \half N\m$ for $N \ge n_1$ on~$A_N$ and, by Cauchy--Schwarz,
\[
    \sum_{l \le \nith}l \sqrt{p_l} \Le \sqrt{N^{1/3} \ex\{\n(1)^2\}}.
\]
Finally,
\[
   \left| 1 - \mu \frac{|J_N| }{M_N} \right| \Eq \frac N{M_N}\,|N^{-1}M_N - \m(1-\z_N)| 
      \Le 2\m^{-1} N^{-7/16}
\]
on~$A_N$, for $N\ge n_1$.

\medskip 
Now, since $G(x)-\e_N \le G_N(x) \le G(x)+\e_N$ for all~$x\ge0$, it also follows for all~$y$
that $G^{-1}(y-\e_N) \le G_N^{-1}(y) \le G^{-1}(y+\e_N)$, and thus that
\bea \label{g-diff}
    |G_N^{-1}(y) - G^{-1}(y)| \le G^{-1}(y+\e_N) - G^{-1}(y-\e_N).
\ena 
Hence it follows that, for any~$\h>0$,
\eqs
  \int_0^{1-\h} |G_N^{-1}(y) - G^{-1}(y)|\,dy
    &\le& \int_0^{1-\h} \{G^{-1}(y+\e_N) - G^{-1}(y-\e_N)\}\,dy \\
    &\le& \int_{1-\h-\e_N}^{1-\h+\e_N} G^{-1}(y)\,dy \Le 2\e_N G^{-1}(1-\h+\e_N) \Eq \ps_N^2.
\ens
Taking $\h := 2\e_N$, this shows that, for $U$ uniformly distributed on $[0,1]$,
\[
    \ex\{|G_N^{-1}(U) - G^{-1}(U)| I[U \le 1-2\e_N] | A_N^* \} \Le \ps_N^2,
\]
and Markov's inequality completes the proof.  Note that, 
since~$G$ is assumed to have finite second moment, $x^2(1-G(x)) = o(1)$ as $x\to\infty$, 
implying that $\e_N G^{-1}(1-\e_N) = o(\e_N^{1/2})$ as $N\to\infty$. 
\ep

\begin{corollary}\label{extremes}
Let~$A_N^*$ be as in Lemma~\ref{G_N-G}.  If~$G$ satisfies Assumption~2 with $b<\infty$, 
then, on~$A_N^*$, 
$$
    \sup_{0\le u\le 1} |G_N^{-1}(u) - G^{-1}(u)| \Eq o(1) \ \mbox{ as }N \to \infty;
$$
if~$G$ satisfies Assumption~2 with $b=\infty$, then
$$
   \sup_{u\colon\,G^{-1}(u) \le x_N} |G_N^{-1}(u) - G^{-1}(u)| \Eq o(1) 
            \ \mbox{ as }N \to \infty,
$$
for~$x_N$ such that $N^{-\a}x_N$ is bounded below as $N\to\infty$ for some~$\a > 0$.
\end{corollary}

\proof
For the first part, let the support of~$G$ be $[a,b]$.  Then for any~$\d > 0$, with \eqref{g-diff}, 
\[
    |G_N^{-1}(u) - G^{-1}(u)| \Le 2\d + \frac{2 \e_N}{g_-[a+\d,b-\d]},
\]
where~$\e_N = O(N^{-1/6})$ is as in Lemma~\ref{G_N-G}, and 
$g_-[c,d] := \inf_{c \le x \le d}g(x)$.  So take $\d = \d_N \to 0$ in
such a way that $\e_N = o(g_-[a+\d_N,b-\d_N])$.

 For the second part, for $N$ large enough that $k(x_0+1)^{-\g} > \e_N$, define $x_{N1} > x_0$
such that $(x_{N1}+1)^\g = k/\e_N$, and choose any $x_N \le x_{N1}$.  Then, uniformly for 
all~$u$ such that $a+\d \le G^{-1}(u) \le x_N$, 
\[
    G^{-1}(u+\e_N) - G^{-1}(u) \Le \frac{\e_N}{\min\{g_-[a+\d,x_0],k(x_N+1)^{-\g}\}}
            \Le \frac{\e_N}{g_-[a+\d,x_0]} + \frac{\e_N}{k(x_N+1)^{-\g}}.
\]
So choose $\d_N = \e_N^{1/2}$ and $x_N = (k\d_N/\e_N)^{1/\g} - 1 \le x_{N1}$,
and $\d_N'$ such that $\e_N = o(g_-[a+\d_N',x_0])$; this gives 
\[
    \sup_{u: G^{-1}(u) \le x_N} |G_N^{-1}(u) - G^{-1}(u)| 
        \Le \frac{2\e_N}{g_-[a+\d_N',x_0]} + 2\d_N + \d_N' \ \to\ 0,
\]
and $x_N N^{-1/(12\g)}$ is bounded below as $N\to\infty$.
\ep

We also need to know that paths of a given length cannot contain too many births.

\begin{lemma}\label{path-numbers}
Suppose that $\lim_{\e\to0} G(\e) = 0$.  Then
there exist $t_* > 0$ such that all individuals of 
generation~$n$ in~$\bZ$ are born after time~$nt_*$, except on an event 
of probability at most $2e^{-n}$.
\end{lemma}

\proof 
Let~$\bZ_n$ denote the number of individuals of generation~$n$ in~$\bZ$, starting
with $\bZ_0=1$.  Then $\ex\bZ_n = \m^n$, and so $\pr[\bZ_n > \{e\m\}^n] \le e^{-n}$.
Now the time elapsed up to generation~$n$ along any given line is a sum of~$n$ independent 
$G$--distributed random variables, 
and the probability that fewer than~$n/2$ of these are
greater than a given value~$\e$ is the binomial probability
\[
   \Bi(n,p)[\lceil n/2 \rceil, n] \Le \{1 + p(z_p-1)\}^n z_p^{-\lceil n/2 \rceil} 
           \Le (4p)^{n/2}  \,,
\]
with $z_p := (1-p)/p$ and $p = G(\e)$.
Hence the probability that, up to generation~$n$, any line takes less that time~$\e n/2$ is at most
\[
    e^{-n} + \exp\{n(\log\m + 1) - \half n\log(1/4G(\e))\}.
\]
Taking $\e > 0$ such that $\log(1/4G(\e)) \ge  2(\log\m + 2)$ makes this probability at most~$2e^{-n}$,
and taking $t_* :=   \e/2$ proves the lemma. 
\ep

\subsection{Controlling the ghosts}\label{coincidence}
We now need to control the differences between the epidemic and branching processes;
we need to show that  ghosts play no significant part.  We begin with
the forward branching process~$Z$.   Recalling from~\Ref{B-def} that 
$W(t) := B'(t)e^{-\l t} \to W$ a.s., we write $e_W := \sup_t \ex W(t)  < \infty$.
Label the individuals of~$Z$ independently and uniformly from~$[N]$ in order
of birth epoch until time~$\t_N$; let $L(t)$ denote the number of times that a label has been
used before, creating an initial ghost, and let $L_+(t) \ge L(t)$ denote the number
of initial ghosts and their descendants whose birth times have been determined by time~$t$.
Finally, let $t_N^\a := \a\l^{-1}\log N$, $\a > 0$.
 
\begin{lemma}\label{ghosts-forward}
Under the above assumptions,
\[
   \pr[\{N^{-1/2}L_+(\t_N) \ge N^{-1/4}\} \cap \{W(\t_N) \ge N^{-1/8}\}]
      \Eq O(N^{-1/8}\log N).
\]
\end{lemma}

\proof
For any of the first~$\nitmi$ indices chosen, the probability that it is a
repeat of an index chosen earlier is at most~$N^{-1/2}(\m+1)$.  Hence, for any
$\a > 0$, writing $T :=  t_N^\a$,
\eqs
    \ex\{L_+(\t_N \wedge t_N^\a)\} &\le& (\m+1)N^{-1/2} \ex\Blb \int_0^T \m e_W e^{\l(T-t)}\,B(dt) \Brb, 
\ens
since an individual born at~$t$ has an expected number of descendants at time~$T$ of at
most $e_W e^{\l(T-t)}$, for each of which the expected number of offspring whose births
are still to come is at most~$\m$.  Hence
\eqs
    \ex\{L_+(\t_N \wedge t_N^\a)\} 
       &\le& (\m+1)N^{-1/2} \m e_W e^{\l T}\ex\Blb B(T)e^{-\l T} + \l \int_0^T e^{-\l t}B(t)\,dt \Brb \\
       &\le& (\m+1)N^{-1/2} \m e_W^2 (1 + \l T)e^{\l T}.
\ens
Thus, choosing $\a = (1+\e)/2$, we have
\[
   \pr[\{N^{-1/2}L_+(\t_N) \ge N^{-1/2+\e}\} \cap \{W(\t_N) \ge N^{-\e/2}\}]
      \Eq O(N^{-\e/2}\log N),
\]
since $\t_N \le t_N^{(1+\e)/2}$ when $W(\t_N) \ge N^{-\e/2}$,
and the lemma follows by taking $\e=1/4$.
\ep

\ignore{
We begin by observing that, for any $w>0$,
\eq\label{repeats}
   \ex\{L(t) I[W(t) \le w]\} \Le \half N^{-1}\{ we^{\l t}\}^2,
\en
so that, for any $0 < \a < 1/2$ and for $\e > 0$, we have
\eqa
   \pr[L(t_N^\a) > 0] &\le& \pr[W(t_N^\a) > N^\e] + \ex\{L(t_N^\a) I[W(t_N^\a) \le N^\e]\} \non\\
     &\le& N^{-\e} e_W + \half N^{2\e + 2\a - 1} \Eq O(N^{-\e}), \label{FG1}
\ena
if we take $\e := \third(1-2\a)$.  It also follows from~\Ref{repeats} that
\[
   \ex\{L(t_N^{(1+\e)/2}) I[W(t_N^{(1+\e)/2}) \le N^\e]\} \Le \half N^{3\e},
\]
implying that 
\eq\label{FG2}
   \pr[\{L(t_N^{(1+\e)/2}) > N^{4\e}\}]
     \Le (\half + e_W) N^{-\e}.
\en
Now, if an initial ghost is born after time~$t_N^\a$, its expected number of 
descendants at time~$t_N^{(1+\e)/2}$ is at most $e_W N^{(1+\e)/2 - \a}$; hence
\[
   \ex\{L_+(t_N^{(1+\e)/2})\,I[L(t_N^{(1+\e)/2}) \le N^{4\e}]\,I[L(t_N^\a) =  0] \} 
             \Eq O(N^{4\e + (1+\e)/2 - \a}).
\]
Taking $\a = 7/16$ and thus $\e = 1/24$, this implies that
\[
   \pr[\{N^{-1/2}L_+(t_N^{(1+\e)/2}) > N^{-5/24}\}\cap\{L(t_N^{(1+\e)/2}) \le N^{4\e}\}\cap\{L(t_N^\a) =  0\}]
     \Eq O(N^{-1/24}),
\]
and hence, from \Ref{FG1} and~\Ref{FG2}, that
\[
   \pr[\{N^{-1/2}L_+(t_N^{(1+\e)/2}) > N^{-5/24}\}] \Eq O(N^{-1/24}).
\]
The final step is to observe that, on $\{W(\t_N) \ge N^{-1/48}\}$, 
$\t_N \le t_N^{(1+\e)/2}$, and hence that $L_+(\t_N) \le L_+(t_N^{(1+\e)/2})$.
}

For the backward branching processes $\bZ_N$ and~$\bZ$, the argument is a little
different, because the identities of the individuals (even if not their labels)
are implicitly recognised during the construction of the branching process~$\bZ_N$;
the choice of a particular value from~$G_N$ may well determine the choice of the
individual in~$J_N$ that gave rise to it, and will certainly do so if the
distribution~$G$ is continuous.  Hence, when constructing~$\bZ_N$, an initial
ghost appears when the same birth time~$t_{Nj,l}$ is sampled from the same
individual~$j$ for the second or subsequent time, and individual~$j$ is represented 
more than once (but without creating ghosts) if several distinct elements 
of~$\{t_{Nj,l},\,1\le l\le j\}$ are sampled.  By Lemma~\ref{backward-progeny}, the 
branching process~$\bZ$ has
the same offspring numbers as~$\bZ_N$ up to~$\nfe$ with probability
$1 - O(N^{-1/8})$, and individuals can also be
identified starting from a realization of~$\bZ$, by using the quantile transformation 
to go from a value sampled from~$G$ to the corresponding value from~$G_N$ (with an 
arbitrary rule for distinguishing individuals that give rise to identical birth times).
Thus the ghosts arise during the joint construction; afterwards, labelling is at random
without replacement from~$J_N$ for the {\it distinct\/} individuals in~$\bZ$
up to time~$\nfe$.

As before, we note that $\bW(t) := \bB(t)e^{-\l t} \to \bW$ a.s.\ as $t\to\infty$. We 
can then write 
$e_{\bW} := \sup_t \ex \{\bW(t) \giv \bB(0) = 1\}  < \infty$ (if the process
is started with $\bB(0)=2$, as from $K$ and~$K'$, the supremum is doubled). 
We let~$\bL(t)$ denote the number of initial
ghosts that have arisen by time~$t$, $\bL_+(t) \ge \bL(t)$ the number of 
initial ghosts and their descendants that have arisen by then, and~$\tL\ut(t)$ the number
of individuals represented at least twice by time~$t$.  We also denote by~$\bL^\th(t)$ the
number of marked individuals and their descendants up to time~$t$, if individuals are
marked independently with probability~$\th$.

\begin{lemma}\label{ghosts-backward}
Let $K$ and~$K'$ be independently chosen at random from~$J_N$, and let
$\h_N' := \h_N\log N$, where~$\h_N = o(N^{-1/24})$ is as in Lemma~\ref{G_N-G}. Then, 
conditional on $A_N^*$,  and starting
the branching process~$\bZ$ either from~$K$ or from both of $K$ and~$K'$, we have 
\eqs
  (1)&& \pr[N^{-1/2}\bL_+(t_N(u)) \ge N^{-3/16} \giv \ffptn \cap A_N]  \Eq O(N^{-1/8}\log N); \\
  (2)&& \pr[N^{-1/2}\bL^{\th(N)}(t_N(u)) \ge N^{-3/16}]  \Eq O(\th(N)N^{5/24}\log N),
\ens
uniformly for all $u \le (\log N)/48$. 
Furthermore, there is a set~$A_N^4 \in \ffptn$ with $\pr[(A_N^4)^c] = O(N^{-1/24})$
such that
\eqs
  (3)&& \pr[N^{-1/2}\tL\ut(t_N(u)) \ge N^{-7/24} \giv  \ffptn \cap A_N^4]  \Eq O(N^{-1/24}),\phantom{XXX}
\ens 
uniformly in the same range of~$u$.
\end{lemma}

\proof
The first and second statements of the lemma are proved in much the same way as Lemma~\ref{ghosts-forward}.
For the first, we note that the probability of the $r$-th individual born being an initial ghost
is at most~$(r-1)/M_N$.  Hence, for any $w > 0$ and $N \ge n_1$, 
\eqs
    \lefteqn{\ex\bigl\{\min\{\bL_+(t_N(u), \bt(\bB,we^{\l t_N(u)}))\} \giv \ffptn \cap A_N \bigr\} }\\
          &&\qquad\Le   M_N^{-1}N^{1/2}we^u \ex\Blb \int_0^{t_N(u)} e_W e^{\l(t_N(u)-t)}\,\bB(dt) \Brb, \\
       &&\qquad\Le 2\m^{-1}we^{2u} e_W^2 (1 + u + \half\log N),
\ens
where $\bt(\bB,v) := \inf\{t\colon\,\bB(t) \ge v\}$.
Thus, and from Lemma~\ref{backward-progeny},
\eqs
    \lefteqn{\pr[\{N^{-1/2}\bL_+(t_N(u)) \ge N^{-1/2 + 5\e/4}\} \cap \{\bt(\bB,N^{\e/2}e^{\l t_N(u)}) > t_N(u)\}
          \giv \ffptn \cap A_N ]} \\
       &&\qquad\Eq O(N^{-3\e/4}e^{2u}\log N).\phantom{XXXXXXXXXXXXXXXXXXXXXXX}
\ens
Since also 
\[
  \pr[\bt(\bB,N^{\e/2}e^{\l t_N(u)}) \le t_N(u)] 
     \Eq \pr[\bB(t_N(u)) \ge N^{\e/2}e^{\l t_N(u)}] \le N^{-\e/2}e_{\bW},
\] 
it follows that, for $u \le \eighth\e\log N$,
\[
    \pr[N^{-1/2}\bL_+(t_N(u)) \ge N^{-1/2 + 5\e/4} \giv \ffptn \cap A_N] \Eq O(N^{-\e/2}\log N),
\]
and the first statement follows by taking $\e = 1/4$.

For the second, we have
\eqs
    \ex\{N^{-1/2}\bL^\th(t_N(u))\} 
          &\le&   N^{-1/2}\th \ex\Blb \int_0^{t_N(u)} e_{\bW} e^{\l(t_N(u)-t)}\,\bB(dt) \Brb, \\
       &\le& \th e^{u} e_{\bW}^2 (1 + u + \half\log N) ,
\ens
and the statement follows from Markov's inequality.

For the third, we begin by noting that the choices of individual in~$\bZ_N$ after~$n$ have been
examined are multinomially $\MN\bigl(n;\,\{\n(j)/M_N,\,j\in J_N\}\bigr)$ distributed, so that the
mean number of individuals that have by then been chosen more than once is at most
\eq\label{n-tries}
       \frac{n^2}2 \sjjn \Bl\frac{\n(j)}{M_N}\Br^2 \Le \frac{n^2}2 \sli M_{Nl} \Bl\frac{l}{M_N}\Br^2.
\en
Let $A_N^4 := \{\sli M_{Nl} (l/M_N)^2 \le 2N^{\e-1}\} \in \ffptn$, and 
suppose that~$N\ge n_1$ as in~\Ref{n1-def}.  
Observe that, since $M_N \ge \half N\m$ on~$A_N$, and
\[
     \ex \Blb \sli M_{Nl} \Bl\frac{2l}{N\m}\Br^2 \Brb \Le \sli Np_l \Bl\frac{2l}{N\m}\Br^2
            \Le 4N^{-1}\ex\{\n(1)^2\}\m^{-2},
\]
we have $\pr[(A_N^4)^c] = O(N^{-\e})$ for any $\e < 1/8$.  
Then, using~\Ref{n-tries},
\[
   \ex \{\tL\ut(t) I[\bW(t) \le N^\e] \giv A_N^4\} 
    \Eq \ex \{\tL\ut(t) I[\bB(t) \le N^\e e^{\l t}] \giv A_N^4\} \Le N^{-1+3\e}\,e^{2\l t} \Le N^{4\e}, 
\]
uniformly in $t \le (1/2\l)(1+\e)\log N$.  Hence, and since $\pr[\bW(t) > N^{\e}] \le e_{\bW}N^{-\e}$,
it follows that, for $u \le \half\e\log N$,
\[
    \pr[N^{-1/2}\tL\ut(t_N(u)) \ge N^{5\e-1/2} \giv A_N^4] \Eq O(N^{-\e}),
\]
giving the third assertion if we take $\e = 1/24$.
\ep

 We now use $\hL(G_N,t)$ to denote
the number of individuals in~$\bZ$, together with their descendants, up to time~$t$,
for which the sample taken from~$G$ to determine their birth time is such that the difference
between it and the corresponding value obtained from~$G_N$ by the quantile
transformation exceeds the threshold~$\ps_N$ defined in Lemma~\ref{G_N-G}.  Note that, 
on~$A_N^*$, 
the expected contribution to~$\hL(G_N,t)$ resulting from the offspring of an individual
born at time~$v < t$ is at most $\m\h_N e_{\bW}e^{\l(t-v)}$, where~$\h_N$ is as
in Lemma~\ref{G_N-G}. The proof of Lemma~\ref{ghosts-backward}(2) then yields the
following corollary.

\begin{corollary}\label{deformation}
In the setting of Lemma~\ref{ghosts-backward}, with $\h_N' = \h_N\log N$ and $\h_N$ as
in Lemma~\ref{G_N-G}, we have
\eqs
   \pr[N^{-1/2}\hL(G_N,t_N(u)) \ge (\h_N')^{1/2} \giv \ffptn\cap A_N^*]  \Eq O((\h_N')^{1/2}).
\ens 
\end{corollary}

\subsection{Main theorem}\label{main-theorem}
We now combine our previous results to prove the main result of Section~\ref{basic}.  For any~$t \ge 0$,
let~$\ss_N(t)$ denote the set of individuals in the epidemic that are still susceptible at time~$t$, 
and write $S_N(t) := |\ss_N(t)|$. Then, for independently and randomly chosen $K$ and~$K'$ in~$[N]$, 
$$
    \ex\{N^{-1}S_N(t) \giv \ffptn\} \Eq \frac{1}{N} \sum_{k=1}^N \pr [ k \in \ss_N(t) \giv \ffptn] \Eq  \pr[K \in \ss_N(t) \giv \ffptn] ,
$$
and similarly 
$$
   \var\{N^{-1}S_N(t) \giv \ffptn\} \Eq \pr[\{K,K'\} \subset \ss_N(t) \giv \ffptn]
      - \{\pr[K \in \ss_N(t) \giv \ffptn]\}^2,
$$
and we use these expressions to show that $N^{-1}S_N(t)$ is close to its expectation, and
to give an asymptotic expression for it.

At time~$\t_N$, the epidemic process has generated a collection~$U_N$ of individuals, whose
birth times, the elements of~$V_N(\t_N)$, are determined, but have not yet occurred,
and which have not yet been labelled (so that some of them may turn out to be ghosts);
labels are assigned to them independently and at random from~$[N]$, and ghosts are then
removed, leaving a labelled set $U'_N \subset U_N$.

A randomly chosen individual~$K$ samples an independent copy of the reversed branching
process~$\bZ$, and uses it to determine its susceptibility process, by way of~$\bZ_N$.
For times to infection, as measured in~$\bZ$-time, not exceeding~$t_N(u+h)$, there is 
a corresponding susceptibility set~$\hU_N(u+h)$, consisting of distinct individuals.
The elements of the set~$\hU_N(u+h)$ are now assigned labels, chosen independently but
{\it without replacement\/} from~$J_N$.  Let $\ee_N(u+h)$ denote the set of elements 
of~$\hU_N(u+h)$ that share labels with members of~$U_N$.  Then $E_N(u+h) := |\ee_N(u+h)|$ 
has conditional expectation $|\hU_N(u+h)|\,|U_N|/N$. 
If $E_N(u+h) = 0$, there is no path of infection from~$i_1$ to~$K$ of $\bZ$-length
less than $\t_N + u + h$.  If $E_N(u+h) > 0$, go through the elements of~$\ee_N(u+h)$
in order of increasing $\bZ$-time, and mark all their progeny in~$\hU_N(u+h)$ as ghosts,
since these elements are also represented as members of~$U_N$, and their infection pre-history 
has already been determined in~$\ffptn$.  Let $\ee'_N(u+h) \subset \ee_N(u+h)$ denote
those elements of $\ee_N(u+h)$ that are not marked as ghosts, and write $E'_N(u+h) := |\ee'_N(u+h)|$. 
For any element~$e$ of~$\ee'_N(u+h)$, let $\t_N+v$ denote the birth time of the corresponding
element of~$U'_N$, let $\bs$ denote the birth time in~$\bZ$ of the element 
of~$\ee'_N(u+h)$, and~$\bs_N$ its corresponding birth time in~$\bZ_N$.   Then~$e$ gives
rise to an infection path from $i_1$ to~$K$ of length $\t_N + v + \bs_N$.  If this is
less than or equal to $\t_N + t_N(u)$ for any~$e$, then~$K \notin \ss_N(\t_N + t_N(u))$;
otherwise, $K \in \ss_N(\t_N + t_N(u))$ unless, possibly, there is an infection path
with $v + \bs_N \le  t_N(u)$ but $v + \bs > t_N(u+h)$.  Using these considerations,
we can deduce an approximation for $\pr[K \in \ss_N(\t_N + t_N(u)) \giv \ffptn]$, and
a similar argument, with two reversed branching processes, leads also to a corresponding
approximation to $\pr[\{K,K'\} \subset \ss_N(\t_N + t_N(u))$.

The proof of the theorem that follows is essentially concerned with quantifying the 
above steps.  In particular, it is to be shown that $\ee_N(u+h) = \ee'_N(u+h)$
with high probability, and that $|\hU_N(u+h)|\,|U_N|/N \sim (\m-1)\bW e^{u+h}$.
Then, for any element~$e$ of~$\ee'_N(u+h)$, we need to show that the corresponding~$v$
is sampled from a distribution close to~$F$, as defined in~\Ref{F_e-def2}, and 
that $t_N(u + h) - \bs_N$ is sampled from
a distribution close to the exponential distribution Exp($\l$) with mean $1/\l$,
in view of~\Ref{expl-dist}.  Assuming that this is the case, it follows that
\eqa
    \pr[v + \bs_N \le t_N(u)] &\sim& e^{-h}\int_0^\infty \l e^{-\l s} F(s)\,ds \non\\
              &=& e^{-h} \frac{\m}{\m-1} \int_0^\infty \l s e^{-\l s}\,G(ds) .  
         \label{mean-integral}                 
\ena
The conditional mean number of such events is therefore asymptotically $\bW e^u m_*$,
where $m_*$ is given in \eqref{m-star-def}.
and a Poisson approximation shows that the probability of none of them occurring is
close to $e^{-\bW e^u m_*}$.  The required approximation to $\pr[K \in \ss_N(\t_N + t_N(u)) \giv \ffptn]$
is then $\ex\{e^{-\bW e^u m_*}\}$.  Finally, the possibility that there is an infection path
with $v + \bs_N \le  t_N(u)$ but $v + \bs > t_N(u+h)$ has to be excluded.

\begin{theorem}\label{main}
Under Assumptions 1 and~2,  there exists an event~$\tA_N \in \ffptn$ such that 
$\pr[\tA_N^c] \to 0$ as $N\to\infty$, for which
\[
   \pr\Bigl[\sup_u |N^{-1}S_N(\t_N + \l^{-1}\{\half\log N + u\}) -\hs(u)| > \e 
                   \Giv \ffptn\cap \tA_N \cap \{\t_N < \infty\}\Bigr] \ \to\ 0
\] 
as $N\to\infty$, where~$\hs$ is the decreasing function given by
\[
   \hs(u) \Def \ex\bigl\{ \exp\{-  \bW e^u m_*\}  \bigr\},
\]
and where~$m_* = \m\l \int_0^\infty s e^{-\l s}\,G(ds)$, as in~\Ref{m-star-def}.
\end{theorem}

\nin{\bf Remark}.\   
It therefore follows that 
$\sup_u |N^{-1}S_N(\l^{-1}\{\log N -\log W + u\}) - \hs(u)| \to_d 0$,
conditionally on $W>0$.  However, in practice,
it may be more reasonable to expect to be able to observe the time~$\t_N$ than it is to know
the value of~$W$, or, equivalently, when the first infection occurred.

\medskip
\proof
By Lemma~\ref{ghosts-forward} and Nerman~(1981), Corollary 5.6 with $\phi_1(t) = \xi(t, \infty)$ 
and $\phi_2(t) =1$, and using the fact that $N^{-1/2}B'(\tau_N) \to 1$ a.s.\ as $N\to\infty$,
we obtain that  
$N^{-1/2}|U_N| \to (\m-1)$ a.s.\ as $N\to\infty$ on $\{W>0\}$; Lemma~\ref{ghosts-forward} 
shows that excluding ghosts has negligible effect on the branching asymptotics.
Thus we can define a set
\eq\label{A_N^5-def}
    A_N^5 \Def \{|N^{-1/2}|U_N| - (\m-1)\,| \le \h_1(N)\}\ \in\ \ffptn,
\en 
where $\h_1(N) \to 0$ and $\pr[(A^5_N)^c] \to 0$ as $N\to\infty$.  Let
\[
   \tA_N \Def A_N^* \cap A_N^4 \cap A_N^5 \cap \{W(\t_N) \ge N^{-1/8}\}.
\]

We wish first to show that, for any $u\in\re$, 
$$
   \pr[K \in \ss(\t_N + t_N(u)) \giv \ffptn \cap \tA_N] \ \sim\ \hs(u),
$$
where~$\hs$ is as stated in the theorem.  To do so we proceed as outlined above.
On~$\tA_N$, we have $|U_N| \sim N^{1/2}(\m-1)$, in view of~\Ref{A_N^5-def}. Then,
by~\Ref{what-cvgce} and Lemma~\ref{ghosts-backward}(1,3), $|\hU_N(u+h)| \sim N^{1/2} e^{(u+h)} \bW$;
Lemma~\ref{ghosts-backward} shows that excluding ghosts and individuals multiply referenced 
has little effect on the branching asymptotics.  The mean number of individuals in~$\hU_N(u+h)$ that share 
a common index with a member of~$U_N$ is thus asymptotic to
$$
   N^{1/2}(\m-1) . N^{1/2}\bW e^{ (u+h)} / N \Eq \bW(\m-1) e^{(u+h)} .
$$ 
We now show that $\pr[\ee_N(u+h) \neq \ee'_N(u+h)] = O(N^{-3/16})$.  Letting $E_N^D(u+h)$
denote the number of descendants of~$\ee_N(u+h)$, 
it follows from Lemma~\ref{ghosts-backward}(2), by taking 
$\th = \th(N) = N^{-1}|U_N|$ and in view of~\Ref{A_N^5-def}, that 
\[
    \pr[E_N^D(u+h) \ge N^{5/16} \giv \ffptn \cap \tA_N] \Eq O(N^{-1/4}\log N).
\]
The conditional probability that any of them is marked by a label from~$U_N$ is thus  
at most of order $O(N^{-1/2+5/16}+N^{-1/4}\log N)  = O(N^{-3/16})$.

Now, because of the random scheme of assignment of labels, any pair in~$\ee_N(u+h)$ 
is associated 
with a random choice of elements~$v$ of~$V_N(\t_N)$ and~$a$ of~$A(t_N(u+h))$, and the 
empirical distributions of the elements of these sets converge, as observed in 
\Ref{F_e-def1}, \Ref{F_e-def2} and~\Ref{expl-dist}.  
Furthermore, the empirical distribution~$\bF_N^{(u+h)}$ of the birth times 
in~$\bZ_N$ corresponding to the elements of~$A(t_N(u+h))$ also converges to the 
exponential~Exp($\l$) distribution with mean~$1/\l$ if $\bW > 0$.  To see this, we argue as follows.
Recalling~\Ref{expl-dist}, let~$\h_2(t)$ be such that $\lti \h_2(t) = 0$ and that
\eq\label{expl-rate}
  \pr\left[ \sup_{s\ge0}\Blm \frac1{\bB(t)} \sum_{r=1}^{\bB(t)} I[a_t(r) \le s] - (1 - e^{-\l s})\Brm
             > \h_2(t) \Giv \bW > 0 \right] \Le \h_2(t).
\en
Then define
\[
    k \Def \Big\lceil \frac{1+\e}{2\l t_*} \Big\rceil, 
\]
where~$t_*$ is as in Lemma~\ref{path-numbers}.  Observe that, 
in view of Corollary~\ref{deformation} and of Lemma~\ref{path-numbers}, 
\[
    \sup_s|\bF_N^{(u)}(s) - (1-e^{-\l s})| \Le \l\ps_N k \log N  + \h_2(t_N(u)) 
           + N^{1/2}(\h_N')^{1/2}/|A(t_N(u))|,
\]
on~$\{\bW > 0\}$, uniformly in~$u \le \half\e\log N$, except on a set of conditional probability 
at most~$(\h_N')^{1/2} + 2N^{-(1+\e)/2\l t_*} + \h_2(t_N(u))$, and that $\sup_t e^{-\l t}|A(t)| < \infty$.

\medskip

At this point, we also need to exclude the possibility that there is an infection path
with $v + \bs_N \le  t_N(u)$ but $v + \bs > t_N(u+h)$.  Corollary~\ref{deformation} shows
that, on~$A_N^*$, the probability of having a path from~$K$ to~$U_N$ 
containing a sample~$\ttt$ from~$G$
such that $|\ttt - \ttt_N| > \ps_N$, where $\ttt_N := G_N^{-1}(G(\ttt))$, before time
$\l^{-1}(1/2 + \e)\log N$ is small for $\e < 1/24$, and the number of births in a path up to
that time is bounded by $c\log N$ in view of Lemma~\ref{path-numbers}, with high
probability.  Hence there has to be at least one pair
$(\ttt,\ttt_N)$ in the path such that $\ttt - \ttt_N > c'$, for $c' = (1/2c\l)\e$, if
$u < \half\l^{-1}\e\log N$ and $\bs-\bs_N \ge \l^{-1}\e\log N - u$.  But this
cannot be the case, for~$N$ large enough, in view of Corollary~\ref{extremes}.

\medskip
\ignore{
Recalling~\Ref{F_e-def1},
let~$\h_1(t)$ be such that $\lti \h_1(t) = 0$ and that
\eq\label{F_e-rate}
  \pr\bigl[ \sup_{s\ge0}\bigl|\,(|V'(t) \cap (0,s]| / |V'(t)|) - F_e(s)\bigr| > \h_1(t) \Giv W > 0\bigr] 
     \Le \h_1(t);
\en
} 

Hence, on the event~$\tA_N$, and conditional on~$\ff_N(u) := \s(\bZ(t),\,0\le t\le t_N(u)) \bigvee \ffptn$, 
the mean number of pairs with common index,
one from~$U_N$ and one from~$\hU_N(u+h)$, that are not ghosts and give rise to an infection
path between $i_1$ and~$K$ of length at most $\t_N + t_N(u)$, is given as in \Ref{mean-integral}
and~\Ref{m-star-def} by
\eqs
    m_N(u,\bW) &\sim& \bW e^{u} m_*;
\ens
of course, the asymptotics are valid also when $\bW = 0$.
Let~$I_{Nj}(u)$ denote the indicator of the event that the label of the $j$-th element of~$U_N$ is
matched with one of the labels assigned to~$\hU_N(u)$, $1\le j\le |U_N|$.  Then,
conditional on~$\ff_N\uu$, 
$(I_{Nj}(u),\,1\le j\le|U_N|)$ is a collection of independent indicator random variables,
each with probability $p_N(u) := |U_N|^{-1}m_N(u,\bW)$; hence it follows by
\BHJ~(1992, (1.23)) that
\eq\label{poisson-1}
   \left| \pr\Bigl[ \sum_{j=1}^{|U_N|} I_{Nj}(u)  = 0 \Giv \ff_N\uu \cap \tA_N \Bigr] 
       -  \exp\{-m_N(u,\bW)\} \right|  \Le p_N(u).
\en  
Thus we deduce that 
\eqs
   \lefteqn{\pr[K \in \ss_N(\t_N+t_N(u)) \giv \ffptn \cap \tA_N ]} \\ 
          &&\ \sim\ \ex\bigl\{ \exp\{-m_N(u,\bW)\} \giv \ffptn \cap \tA_N \bigr\} 
           \ \sim\ \ex\bigl\{ \exp\{-  \bW e^{ u} m_*\}  \bigr\} \ = \ \hs(u).
\ens
But this means that
\eqa
     \hs_{N}(u) &:=& \ex\bigl\{N^{-1}S_N(\t_N + \l^{-1}\{\half\log N + u\}) 
                 \giv \ffptn \cap \tA_N \bigr\} \non\\
           &=&  \pr[K \in \ss_N(\t_N+t_N(u)) \giv \ffptn \cap \tA_N ]  \ \sim\ \hs(u) \label{y-mean}
\ena
also.      

The argument for approximating the probability that both $K$ and~$K'$ belong to~$\ss_N(\t_N+t_N(u))$
runs in much the same way.  The limiting random variable for the branching process~$\bZ$
starting with two individuals can be expressed as~$\bW_1+\bW_2$, where the two are
independent copies of~$\bW$, and the sizes of the corresponding sets $\hU_N\ui(u)$ and~$\hU_N\ut(u)$
are asymptotically $N^{1/2}\bW_1 e^{ u}$ and~$N^{1/2}\bW_2 e^{ u}$ respectively.  
 We write $\tI_{Nj}(u) = (1,0)$ if the $j$-th element of~$U_N$ is matched with a label
 associated with~$\hU_N\ui(u)$, and $(0,0)$ otherwise; similarly, $\tI_{Nj}(u) = (0,1)$ if matched with a label
 associated with~$\hU_N\ut(u)$ and $(0,0)$ otherwise.    Then both $K$ and~$K'$ belong 
to~$\ss_N(\t_N+t_N(u))$ if $\sum_{j=1}^{|U_N|} \tI_{Nj}(u) = (0,0)$.
The multivariate analogue of the Poisson approximation~\Ref{poisson-1}
(Roos, 1999, Theorem~1) gives  
\eqa
   \lefteqn{\left| \pr\Bigl[ \sum_{j=1}^{|U_N|} \tI_{Nj}(u) \in\nat^2 \Giv \ff_N\uu \cap \tA_N \Bigr] 
          - \exp\{-m_N(u,\bW_1)-m_N(u,\bW_2)\} \right| } \non\\
         && \Le c|U_N|^{-1}\{m_N(u,\bW_1)+m_N(u,\bW_2)\}, \phantom{XXXXXXXXXXXXX}
              \label{poisson-2}
\ena 
for a universal constant~$c$. Hence, as before, 
\eqa
   \lefteqn{\pr[\{K,K'\} \subset \ss_N(\t_N+t_N(u)) \giv \ffptn \cap \tA_N ] }\non\\
           &&\ \sim \ex\bigl\{ \exp\{-m_N(u,\bW_1)-m_N(u,\bW_2)\} 
                                                        \giv \ffptn \cap \tA_N \bigr\} 
           \ \sim\ \{\hs(u)\}^2, \label{joint-prob}
\ena
by the independence of $\bW_1$ and~$\bW_2$.  But the joint probability can also
be written as
\[
     \ex\bigl\{\Bl N^{-1}S_N(\t_N + \l^{-1}\{\half\log N + u\}) \Br^2
                 \giv \ffptn \cap \tA_N \bigr\},  
\]
so that it follows from \Ref{y-mean} and~\Ref{joint-prob} that
\eq\label{y-var}
    \var\bigl\{N^{-1}S_N(\t_N + \l^{-1}\{\half\log N + u\}) 
                 \giv \ffptn \cap \tA_N \bigr\} \ \sim\ 0.
\en
It now follows, by a standard argument, that, for any $\e>0$, conditional on $\ffptn\cap \tA_N$,
\[
   \pr\Bigl[\sup_u |N^{-1}S_N(\t_N + \l^{-1}\{\half\log N + u\}) -\hs(u)| > \e \Giv \ffptn\cap \tA_N\Bigr] \ \to\ 0
\]  
as $N\to\infty$, and the theorem follows.
\ep

Because of the factor~$\l^{-1}$ in the definition of~$t_N(u)$, the quantity~$\hs(\l t)$ should
match the solution~$s(t)$ of~\Ref{K-McK}. To see that this is so, note that, by considering
the possibilities for the offspring of the first individual in~$\bZ$, 
$\ps(\th) := \ex\{e^{-\th\bW}\}$ satisfies the equation
\eq\label{branching-eqn}
   \ps(\th) \Eq \exp\Blb -\m\int_0^\infty (1 - \ps(\th e^{-\l w}))\, G(dw) \Brb.
\en
Substituting $\th = m_* e^{\l t}$, writing
\eq\label{s-to-psi}
    s(t) \Eq \hs(\l t) \Eq \ps(m_* e^{\l t})
\en 
and taking logarithms recovers equation~\Ref{deterministic}, 
 with $\m G(du)$ in place of~$\b(v)\,dv$.  
As for~\Ref{deterministic}, equation~\Ref{branching-eqn}
has many solutions, since, if~$\ps(\th)$ is a solution, so is $\ps_\a(\th) := \ps(\a\th)$, 
for any fixed~$\a > 0$. The condition $\ps(0) = 1$, equivalent to $s(-\infty) = 1$, is
satisfied by all~$\ps_\a$.  The relevant choice of solution
to~\Ref{branching-eqn} is determined by matching $\ex\bW$ with~$-\ps'(0)$, or, in terms
of~\Ref{deterministic}, with $(m_* \l)^{-1}\lim_{t\to -\infty}e^{-\l t}(-Ds(t))$.  A renewal 
equation for $\ex\{\bB(t)e^{-\l t}\}$ gives the solution as
\eqs
    \ex\bW &=& \lti \ex\{\bB(t)e^{-\l t}\} \Eq \Blb \l\m \int_0^\infty ve^{-\l v}\,G(dv) \Brb^{-1}
       \Eq \frac1{m_*},
\ens
by the key renewal theorem.  Thus Theorem~\ref{main} can be interpreted as a formal
justification of the stochastic basis for the Kermack--McKendrick epidemic as described in 
Metz~(1978), Section~4, under 
 assumptions that are slightly more general, in that the point processes~$\xi$ are not
required to be doubly stochastic, but are in some respects more restrictive as regards the 
choice of~$\b$.
Since~$\ps$ is identified as the Laplace transform of a probability
distribution, it is an analytic function in $\Re(\th) > 0$, which, with~\Ref{s-to-psi}, 
proves Conjecture~(f) in Metz~(1978), p.120. 

\section{Refinements}\label{refinements}
 \setcounter{equation}{0}
 \subsection{Multitype epidemics}\label{multitype}
Very similar arguments can be carried through for epidemics in populations consisting
of individuals of more than one type.  Suppose that there are a finite number~$d$ 
of different types,
with~$N_l$ individuals of type $l$, $1\le l\le d$, where $N_l \in \{ \lfloor N\p_l \rfloor, 
\lceil N\p_l \rceil\}$, $\slid N_l = N$ and $\slid \p_l = 1$.  
 Assume that type~$l$ individuals have independent and identically
distributed point processes $\xi\ul_i$, $1\le i\le  N_l$ , on $[d]\times\re_+$,
with mean measures
\eq\label{mv-intensity}
    \ex\{\xi\ul_1(k,du)\} \Eq \m_{lk} G_{lk}(du), 
\en
where $\int_0^\infty G_{lk}(du) = 1$ for all $1\le k,l\le d$.   Then an epidemic process can
be constructed in the population, just as in the single type case, by beginning with a
multitype branching process constructed from independent realizations of the~$\xi\ul_1$,
$1\le l\le d$, and then using random labelling within the members of each type
to determine which transitions are to be retained in the epidemic process.
The approximation arguments are very much as before. 
Asymptotically exponential growth and the analogues of \Ref{F_e-def1} and~\Ref{expl-dist},
together with an asymptotically stable type distribution,
hold in~$L_1$ in the multitype setting.  The asymptotic
statements that we use in this section are all justified by Theorem~7.3 of Jagers~(1989), 
who proves~$L_1$ approximation for a wide variety of characteristics of the
branching process in an even more general setting.

\medskip
\nin{\bf Remark.}\
It is perhaps more natural, especially when comparing the spread of the same epidemic
in populations with different compositions of types, to assume a fixed value for
the measures $\a_{lk}(du) := \m_{lk}G_{lk}(du)/\p_k$, rather than supposing that~$\m_{lk}G_{lk}$
remains the same for all~$N$.  The quantity~$\a_{lk}(du)$ can be interpreted as representing 
the infection intensity
measure of contacts with type~$k$ individuals made by a type~$l$ individual, in a
population consisting entirely of individuals of type~$k$.  At least in Poisson
process contact models, this would suggest taking  
$\ex\{\xi^{(l,N)}_1(k,du)\} = \a_{lk}(du)N_k/N$ in a population of the composition
given above, implying that $G_{lk}(du) = \a_{lk}(du)/\a_{lk}(\re_+)$ is fixed
for all~$N$, but that $\m_{lk}^{(N)} = \a_{lk}(\re_+) N_k/N$ may vary with~$N$.
This differs from~\Ref{mv-intensity} inasmuch as $N_k/N$ is not exactly equal
to~$\p_k$.  As in the single--type model, this minor
difference entails no change in the theorems that we prove. 

\medskip
We now assume that the matrix~$\m$ is irreducible, and that 
the distribution functions~$G_{lk}$ all satisfy Assumption~2;  suppose also that the largest eigenvalue
of~$\m$ is larger than~$1$, and write
\[
   \m_{lk}(s) \Def \m_{lk} \int_0^\infty e^{-su}\,G_{lk}(du).
\]
Then the branching process has as Malthusian parameter the value~$\l > 0$ for
which $\m(\l)$ has largest eigenvalue~$1$. We write $\z^T$ and~$\h$ for the
positive left and right eigenvectors of~$\m(\l)$ associated with eigenvalue~$1$,
normalized such that $\z^T 1 = \z^T\h = 1$.    Let $B'(t) := (B_l'(t),\,
1\le l\le d)$ denote the numbers of individuals of each type born up to time~$t$.
Then, if the branching process starts from a single individual of type~$i$,
\eq\label{B-vec-cvgce}
     B'(t) e^{-\l t} \ \to\ W\uii\z \ \mbox{in}\ L_1
\en
as $t\to\infty$. Here, $W\uii$ is a random variable whose Laplace transform
$\ps\uii(s) := \ex\{e^{-s W\uii}\}$ satisfies the implicit equations
\eq\label{LT-vec-forward}
    \ps\ul(s) \Eq \ex\Blb \exp\Bl\skid \int_0^\infty \log\ps\uk(se^{-\l v})\xi\ul(k,dv) \Br \Brb,
     \qquad 1\le l\le d,
\en
with $\ex W\uii = \h_i/m_*\ui$ and
\eq\label{m-star-1-def}
   m_*\ui \Def \l\z^T(-D\m(\l))\h; 
\en
note that
\[
    (-D\m(\l))_{lk} \Eq \m_{lk} \int_0^\infty ue^{-\l u}\,G_{lk}(du),
\]
and that $m_*\ui/\l$ is the multitype mean age at child bearing (Jagers~(1989), p.195).
Letting $V'_l(t)$ denote the set of times until birth of the unborn type~$l$ offspring
of individuals born before~$t$, it follows also that 
\eq\label{unborn}
   e^{-\l t}|V'_l(t)| \ \to\ W\uii c_l\ \mbox{in}\ L_1,
\en
with
\eqa
    c_l &:=& \skid \z_k  \m_{kl}\int_0^\infty (1 - e^{-\l v})\,G_{kl}(dv) \non\\
        &=& \skid \z_k (\m_{kl} - \m_{kl}(\l)) \Eq \skid \z_k \m_{kl} - \z_l,
     \label{c_l-def}
\ena
and that, on $W\uii > 0$, 
\eq\label{F_l-def1}
    \ex\uii\!\Bl\sup_s\Bigl| |V_l'(t)\cap (s,\infty)|/|V_l'(t)| - (1 - F_l(s)) \Bigr|\Br \ \to\ 0,
\en
where
\eq\label{F_l-def2}
    1 - F_l(s) \Def  c_l^{-1}\skid \z_k \m_{kl} \int_s^\infty (1 - e^{-\l(v-s)})\, G_{kl}(dv)\,,
\en
replacing \Ref{F_e-def1} and~\Ref{F_e-def2}.

The backward branching process is similar, but has Poisson point processes $\bxi\ul$
with intensity $\m_{kl}G_{kl}(du)$ at $(k,u) \in [d]\times\re_+$.  
The matrix~$\bm(s)$ is given by~$\m(s)^T$, so that the Malthusian
parameter is still~$\l$, but the left and right eigenvectors at~$\l$ are swapped;
the normalized versions are $\hz^T := \h^T/H$ and $\heta := H\z$, where
$H := \skid \h_k$. The backward random variables
$\bW\ul := \lim_{t\to\infty}e^{-\l t}\skid\bB_k(t)$ corresponding to the initial
conditions $1\le l\le d$ now have means $\heta_l/m_*\ui = H\z_i/m_*\ui$,
and their Laplace transforms~$\baf\ul$ satisfy the equations   
\eq\label{LT-vec-backward}
  \baf\ul(s) \Eq  \exp\Bl -\skid \m_{kl}\int_0^\infty(1 - \baf\uk(se^{-\l v}))\,G_{kl}(dv) \Br,
    \qquad 1\le l\le d.
\en
As in~\Ref{expl-dist},
the empirical distribution of the ages at time~$t$ of $l$-individuals born before~$t$ 
also converges in~$L_1$ to Exp($\l$). 

Now suppose that the forward branching process starts with a single type~$i$ individual.
Define $\t_N := \inf\{t > 0\colon\, \slid B'_l(t) \ge \nit\}$, so that $W\uii e^{\t_N} \sim \sqrt N$
as $N\to\infty$, from~\Ref{B-vec-cvgce}, and $|V'_l(\t_N)| \sim c_l\sqrt N$, $1\le l\le d$, 
from~\Ref{unborn}. Then run the backward branching process starting with a single type~$i'$
individual; at time $t_N(u) := \l^{-1}(\half\log N + u)$, as in~\Ref{t_N-def}, we have
$\bB(t_N(u)) \sim \sqrt N \bW\uid e^{\l u}\hz$.  Hence the mean number of pairs 
consisting of one element~$v$ of~$V'_l(\t_N)$ and one type~$l$ individual~$w$ born before~$t_N(u)$
in the backward branching process, such that $v$ is less than the age of~$w$ at~$t_N(u)$, is
asymptotically given by
\eqs
    \lefteqn{\{c_l\sqrt N\}\,\{ \sqrt N \bW\uid e^{\l u}\hz_l\}\, \int_0^\infty \l e^{-\l s} F_l(s)\,ds }\\
       &&\qquad \Eq N \bW\uid e^{\l u} \hz_l\, \int_0^\infty \l v e^{-\l v}\skid \m_{kl} \z_k\, G_{kl}(dv).
\ens
Thus, when the individuals corresponding to the~$V'_l(\t_N)$ and the type~$l$ individuals in
the backward branching process are randomly labelled in constructing the epidemic process,
the mean number of such pairs that have the same labels is asymptotically given by
\[
   \p_l^{-1}\bW\uid e^{\l u} \hz_l\, \int_0^\infty \l v e^{-\l v}\skid \m_{kl} \z_k\, G_{kl}(dv),
\]
and hence the probability that there is no such pair of any type~$l$, $1\le l\le d$, is
asymptotically given by $\exp\{-\bW\uid e^{\l u} m_*\ut\}$, where
\eq\label{m-star-2-def}
   m_*\ut \Def \int_0^\infty \l v e^{-\l v}\skid\slid  \z_k\, \a_{kl}(dv) \h_l/H.
\en
Arguing as in the case of a single type, we have the following theorem, in which $\ffptn$ denotes the
precise analogue of the $\s$-algebra having the same name in the single type case, and
$S_{Nl}(t)$ is the number of type~$l$ susceptibles at time~$t$.

\begin{theorem}\label{main-multitype}
Suppose that the multitype forward branching process is supercritical and has offspring
distributions with finite second moments; suppose also that Assumption~$2$ holds for each~$G_{lk}$. 
Then there exists an event~$\tA_N \in \ffptn$ such that $\pr[\tA_N^c] \to 0$ as $N\to\infty$,
for which
\[
   \pr\Bigl[\sup_u |(Np_l)^{-1}S_{Nl}(\t_N + \l^{-1}\{\half\log N + u\}) -\hs_{l}(u)| > \e 
                   \Giv \ffptn\cap \tA_N \cap \{\t_N < \infty\} \Bigr] \ \to\ 0
\] 
as $N\to\infty$, where~$\hs_{l}$ is the decreasing function given by
\[
   \hs_{l}(u) \Def \baf\ul( e^u m_*\ut)  ,
\]
where the $\baf\ul$ satisfy~\Ref{LT-vec-backward} with $-D\baf\ul(0) = H\z_l/m_*\ui$,
and where~$m_*\ui$ is defined in~\Ref{m-star-1-def} and~$m_*\ut$ in~\Ref{m-star-2-def}.
\end{theorem}

\subsection{A configuration model}\label{configuration}
In this section, we consider a different model of epidemic spread.  In those considered 
so far, an infected individual chooses to infect a number of randomly chosen individuals,
and the individuals chosen are not taken into account in this choice.
Now we suppose that pairs of individuals are either acquainted with one another 
or are not, so that acquaintanceship determines a graph on the set of individuals,
and we assume that infectious contacts can only be made between graph neighbours.
This yields a more symmetric description of the contact process, and, as a result,
the forward and backward branching approximations can be expected to look more similar. 
We shall, for simplicity, assume that there is a finite, $N$-independent upper bound~$K$
on the number of acquaintances that an individual may have; note that this immediately
rules out any Poisson distribution of offspring in an approximating branching
process, so that the backward branching processes from such a model have to be
different from those in the previous sections.  

To make further progress, we assume
that the acquaintanceship graph is nonetheless rather randomly constituted within
the population, according to the following construction. 
We assume that~$N_k$ members of the population are
`type~$k$' individuals, who have exactly~$k$ acquaintances, with $\skK N_k = N$ and 
$N_k \in \{ \lfloor N\p_k \rfloor, \lceil N\p_k \rceil\}$, for fixed $\p_1,\ldots,\p_K$,
and with $M := \skK kN_k$ even.  Think of a type~$k$ individual as
having~$k$ half--edges, and join the half--edges into edges by means of a random
matching of the~$M$ half--edges, determining the acquaintanceship graph.  This
graph may have some loops and multiple edges, but they are few, and we shall
ignore their effects. Thus the method of assigning which
individuals are acquaintances remains essentially random, but the
propensities of each individual are respected when determining whether
they are acquainted or not.
We then assume that an infected type~$k$ individual makes contact
with a given type~$l$ acquaintance at a random time after infection that has 
(possibly defective) distribution
function~$G_{kl}$ and is independent of all other contact times; we suppose also that
a type~$k$ individual remains infectious for a random time with (possibly defective) 
distribution~$\Phi_k$,
again independently of everything else.  If we specialize to the case where the
distributions~$G_{kl}$ are all identical and equal to Exp($\a$), and that the~$\Phi_k$
are all identical and equal to Exp($\b$), then the model of Volz~(2008) (in the
case of a finite number~$K$ of possible contact numbers) is recovered. 

   As in the previous models, the key effort lies in determining the probability that
an initially chosen individual infects another randomly chosen individual before
a specified time~$t$.  To do so, construct the association graph by starting from the
initial individual as root vertex, and matching its half--edges by random choice from the 
set of all half--edges; then attach the infectious period to the initial individual, and
the lengths of time to potentially infectious contact to the edges.  This yields a set
of infected vertices, together with the times of their infection, some of which may be infinite.
Now continue by matching the remaining half--edges associated with the first of these
vertices (if any) to be infected, attaching the infectious period to the chosen vertex,
and adding the lengths of time to potentially infectious contact (infinite, if longer
than the infectious period) for each edge to the time of infection of the chosen vertex, 
so as to yield the times of infection of newly infected vertices; this augments the set 
of infected vertices. Proceed in this way, always choosing for development the infected 
vertex with unmatched half--edges that has the smallest time of infection, until
the first time that either at least~$\nit$ vertices have been infected or the infection
dies out. In the former case, there remains a set of infected vertices whose subsequent
contact history has not been explored.
If a half--edge is picked for a second or subsequent time, ignore the choice and
re-sample until a new one is chosen; if a vertex is chosen that has already been infected,
ignore it for future development.  As in the previous arguments, for the lengths
of time in which we are interested, there are a few such repeated samples, but few
enough that they can be ignored.  

For the susceptibility graph seen backwards from a randomly chosen individual, carry out 
essentially the same procedure for a specified time; the only difference is the vertex 
to which the infectious period is attached, being that of the child, rather than the
parent. Half--edges that have previously been used, including those that were used
in the forward process, are discarded and re-sampled; the half--edges that are 
associated with the set of infected but unexplored vertices from the forward phase are 
still available for choice, and are those that close chains of infection.

If repeats are ignored, the infection process as
seen from the initial individual becomes a branching process with~$K$ types.  
In the branching process, a type~$k$ individual (other than the initial individual) 
has~$k-1$ offspring, corresponding to the~$k-1$ half--edges that remain to be connected
after a type~$k$ individual has been encountered in growing the association graph, 
and each of these is of type~$l$ with probability $lp_l/m$, where $m = \sum_{l'=1}^K
l'p_{l'}$, chosen from the size--biased transform of the frequency distribution
$(p_l',\,1\le l'\le K)$.  As before, the difference between the process with this distribution and 
that with offspring probabilities $lN_l/M$ is negligible for our purposes.  The type~$k$
individual also has an infectious period randomly assigned to it from the distribution~$\Phi_k$,
and the times to contact along the different edges are assigned independently from the
appropriate distributions~$G_{kl}$.  This yields an age--dependent multi-type branching 
process, in which times to birth may be infinite (if the sampled time to contact is
itself infinite, or exceeds the infectious period of the parent), and the times
of birth of the descendants of a given individual are dependent, because they are
finite only if they do not exceed the infectious period of the common parent.

Seen from the randomly chosen individual, the backward branching process is very 
much the same.  The offspring distribution is identical, but the infection times of
the offspring of a given individual, although having the same marginal distributions 
as before, are now independent, because the relevant infectious period, determining
whether a contact results in infection, is that of the 
child, and not of the parent.  Because the basis of the construction is the fixed
set of half--edges, the problems that arose in Section~\ref{one-type}, because 
the offspring distribution of the backward branching process was not fixed for all~$N$, 
no longer appear (except for the trivial differences between $lp_l/m$ and~$lN_l/M$);
more importantly, choosing the contact times for type~$k$ -- type~$l$ contacts
{\it independently\/} from~$G_{kl}$ and the infectious periods independently from
the~$\Phi_k$ means that the times to birth in the backward branching process have
distributions that do not depend on~$N$, so that there is no need for an analogue
of Corollary~\ref{extremes}, and hence no special assumptions about the tails of
the~$G_{kl}$ need to be made.  Of course, the offspring distributions of the
different types are bounded, so that the corresponding moment conditions are 
automatically satisfied.

The argument now proceeds much as for the multitype process of the previous
section.  Once again, the asymptotic statements for the branching processes are
justified by Jagers~(1989), Theorem~7.3. The matrix~$\m$ is defined analogously by
\[
   \m_{lk}(s) \Def (l-1)\{kp_k/m\}\,\int_0^\infty e^{-su}\,(1 - \Phi_l(u))G_{lk}(du)
     \ =:\ (l-1)\{kp_k/m\} U_{lk}(s),
\]
say, and we write $\m_{lk} := \m_{lk}(0)$;  note that~$\m_{lk}$ need 
no longer be the expected number of offspring, since~$U_{lk}(0)$ is typically
less than~$1$.  Because of the factor~$(l-1)$, $\m(s)$ is reducible.
Supposing that all the~$p_k$ and all the~$U_{lk}(s)$ are positive, we
can write the irreducible non-negative matrix~$\m\ui(s)$, obtained 
from~$\m(s)$ by removing the first row and column, as $D_1 U\ui(s) D_2$, 
where $D_1 := \diag(1,2,\ldots,K-1)$ and $D_2 := m^{-1}\diag(2p_2,\ldots,Kp_K)$.
Assume that the matrix~$\m\ui(0)$ has dominant eigenvalue larger
than~$1$, and define the Malthusian parameter~$\l$ to be such 
that~$\m\ui(\l)$ has dominant eigenvalue equal to~$1$; let $\z\uit$ and~$\h\ui$
be associated left and right eigenvectors.
Then the left and right eigenvectors of~$\m\ui(\l)$ with eigenvalue~$1$ are
given by $\z^T := Z^{-1}(\z\uit\m(\l)\e\ui,\z\uit)$ and
$\h := H^{-1}(0,\h\uit)^T$, where $Z$ and~$H$ are chosen so that $\z^T 1 = \z^T\h = 1$;
here, $\e\ui$ denotes the first coordinate vector.
   
Let $B'(t) := (B_l'(t),\, 1\le l\le d)$ denote the numbers of individuals of 
each type born up to time~$t$; then
\eq\label{B-vec-cvgce-C}
     B'(t) e^{-\l t} \ \to\ \tW\uii\z \ \mbox{in}\ L_1
\en
as $t\to\infty$, if the initial individual has type~$i$. The distribution of~$\tW\uii$ is not quite
the one that would be expected when starting the branching process with a typical type~$i$
individual, because the {\it initial\/} type~$i$ individual has~$i$ offspring, instead of~$i-1$.
However, it can easily be deduced from the Laplace transforms $(\ps\ul(s),\,1\le l\le K)$ of
the limiting random variables for the branching process that has all individuals, including the
initial one, obeying the same rules.  These solve a system of implicit equations that can be
deduced from~\Ref{LT-vec-forward}.  Here, the quantity within the expectation
in~\Ref{LT-vec-forward} can be written as
\[
    \prod_{r=1}^{l-1} \Blb \Bigl(\ps^{(K_r)}(se^{-\l V_r})\Bigr)^{I[V_r \le T]} \Brb,
\]
where~$T$ denotes the infectious period of the type~$l$ individual,
and $K_r$ denotes the type and~$V_r$ the contact time of the $r$-th of his $(l-1)$ acquaintances.
$T$ and $(K_r,\,1\le r\le l-1)$ are independent, and, given $K_r=k$, $V_r$ is drawn independently
of everything else from the distribution~$G_{lk}$.  Thus~\Ref{LT-vec-forward} reduces here to 
the system 
\eq\label{LT-vec-forward-C}
    \ps\ul(s) \Eq \int_{[0,\infty]}  \Blb 
      m^{-1}\skK kp_k\Bl  \int_{[0,t]} \ps\uk(se^{-\l v})G_{lk}(dv) + [1-G_{lk}(t)] \Br \Brb^{l-1} \,\Phi_l(dt),
\en
for $1\le l\le d$, with $-( D\ps\ul)(0) = \h_l/m_*\ui$ and
\eq\label{m-star-1-def-C}
   m_*\ui \Def \l\z^T(-D\m)(\l)\h; 
\en
the Laplace transform of the distribution of~$\tW\ul$ is then given by
\eq\label{LT-vec-forward-C-start}
    \ex\Bigl\{e^{-s\tW\ul}\Bigr\} \Eq \int_{[0,\infty]}  \Blb 
      m^{-1}\skK kp_k\Bl  \int_{[0,t]} \ps\uk(se^{-\l v})G_{lk}(dv) + [1-G_{lk}(t)] \Br \Brb^l \,\Phi_l(dt),     
\en
for $1\le l\le d$.  Here, we have
\[
    (-D\m)(\l)_{lk} \Eq 
            (l-1)\{kp_k/m\}\,\int_0^\infty ue^{-\l u}\,(1 - \Phi_l(u))G_{lk}(du). 
\]
Letting $V'_l(t)$ denote the set of times until birth of the unborn type~$l$ offspring
of individuals born before~$t$, it follows also that, 
if the initial individual is of type~$i$, then
\eq\label{unborn-C}
   e^{-\l t}|V'_l(t)| \ \to\ \tW\uii c_l \quad \mbox{in}\ L_1,
\en
with
\eqa
    c_l &:=& \skK \z_k  (k-1)\{lp_l/m\}\int_0^\infty (1 - e^{-\l v})\,(1-\Phi_k(v))G_{kl}(dv) \non\\
        &=& \skK \z_k (\m_{kl}(0) - \m_{kl}(\l)) \Eq \skK \z_k \m_{kl} - \z_l.
     \label{c_l-def-C}
\ena
Furthermore, on $\tW\uii > 0$,  
as for \Ref{F_l-def1} and~\Ref{F_l-def2},
\eq\label{F_l-def1-C}
   \ex\uii\!\Bl \sup_s\Bigl| |V_l'(t)\cap (s,\infty)|/|V_l'(t)| - (1 - F_l(s)) \Bigr|\Br \ \to\ 0,
\en
where
\eq\label{F_l-def2-C}
    1 - F_l(s) \Def  
     c_l^{-1}\skK \z_k (k-1)\{lp_l/m\} \int_s^\infty (1 - e^{-\l(v-s)})\, (1-\Phi_k(v))G_{kl}(dv)\,.
\en

The backward branching process is similar; we now have 
\[
  \hmu_{lk}(s) \Def (l-1)\{kp_k/m\} U_{kl}(s),
\] 
once again reducible, with $\hmu\ui(s) = D_1 U\uit(s) D_2$ irreducible.
It can be checked that
the Malthusian parameter is still~$\l$. The matrix $\hmu\ui(\l)$ has left and right
eigenvectors $\hz\uit = \h\uit D_2 D_1^{-1}$ and $\heta\ui = D_2^{-1}D_1\z\ui$
with eigenvalue~$1$, and the corresponding left and right eigenvectors of~$\hmu(\l)$
are given by $\hz^T = \hZ^{-1}(\h\uit D_2 U^T(\l) D_2 \e\ui,\hz\uit)$
and $\heta = \hH^{-1}(0,\heta\uit)^T$, where $\hZ$ and~$\hH$ are chosen to make
$\hz^T1 = \hz^T\heta = 1$; in particular, it follows that $\hZ\hH = 1$,
and that the value of $m_*\ui$ deduced from~\Ref{m-star-1-def-C} for the backward process
is the same as~$m_*\ui$.
The limiting random variable~$\tbW\uii$ for the backward
process starting with a single individual of type~$i$, satisfying
\eq\label{B-vec-cvgce-C-back}
     \bB'(t) e^{-\l t} \ \to\ \tbW\uii\hz \ \mbox{in}\ L_1\,,
\en
once again has a distribution whose Laplace transform~$\tbaf\uii$ can be found from the solutions
to a set of implicit equations belonging to the backward branching process whose individuals,
including the initial individual, all follow the same rules.  This branching process
has offspring that behave independently of one another as regards both type and time of birth, 
so that, denoting the Laplace transforms of the limit random variables
with the different initial conditions by $(\baf\ul,\,1\le l\le K)$, we have
$\baf\ul(s) \Eq \{\baf_0\ul(s)\}^{l-1}$, where the~$\baf_0\ul$ satisfy 
the equations   
\eqa
  \baf_0\ul(s) &=&  m^{-1}\skK kp_k 
    \Blb \int_{(0,\infty)} \{\baf_0\uk(se^{-\l v})\}^{k-1}\,(1-\Phi_k(v))G_{kl}(dv)  + (1 - U_{kl}(0))
       \Brb \non\\
   &=& 1 - m^{-1}\skK kp_k \int_{(0,\infty)}(1 - \{\baf_0\uk(se^{-\l v})\}^{k-1})(1-\Phi_k(v))G_{kl}(dv),
     \label{LT-vec-backward-C-back}
\ena
for $1\le l\le d$.
Since~$(-D\baf\ul)(0) = \heta_l/m_*\ui$, the side condition for solving~\Ref{LT-vec-backward-C-back}
is $(-D\baf_0\ul)(0) = \heta_l/\{(l-1)m_*\ui\} = \hH^{-1}\{m/lp_l\}\z_l$, $l \ge 2$, with 
$\baf_0\ui(s) = 1$ for all~$s$. 
The Laplace transform~$\tbaf\uii$ of~$\tbW\uii$ is then given by $\{\baf_0\uii\}^i$.
As in~\Ref{expl-dist},
the empirical distribution of the ages at time~$t$ of $l$-individuals born before~$t$ 
also converges to Exp($\l$).  

Now suppose that the forward branching process starts with a single type~$i$ individual
(having~$i$ offspring).
Define $\t_N := \inf\{t > 0\colon\, \sum_{l=1}^K B'_l(t) \ge \nit\}$, so that $\tW\uii e^{\t_N} \sim \sqrt N$
as $N\to\infty$, from~\Ref{B-vec-cvgce-C}, and $|V'_l(\t_N)| \sim c_l\sqrt N$, $1\le l\le K$, 
from~\Ref{unborn-C}. Then run the backward branching process starting with a single type~$i'$
individual; at time $t_N(u) := \l^{-1}(\half\log N + u)$, 
we have $\bB(t_N(u)) \sim \sqrt N \tbW\uid e^{\l u}\hz$.  
Hence the mean number of pairs of
individuals consisting of an element~$v$ of~$V'_l(\t_N)$ and a type~$l$ individual~$w$ born before~$t_N(u)$
in the backward branching process, such that $v$ is less than the age of~$w$ at~$t_N(u)$, is
asymptotically given by
\eqs
    \lefteqn{\{c_l\sqrt N\}\,\{ \sqrt N \tbW\uid e^{\l u}\hz_l\}\, \int_0^\infty \l e^{-\l s} F_l(s)\,ds }\\
       &&\qquad \Eq N \tbW\uid e^{\l u} \hz_l\, \int_0^\infty \l v e^{-\l v}\skK  \z_k(k-1)\{lp_l/m\}\, 
          (1-\Phi_k(v))G_{kl}(dv).
\ens
Any such pair is realized as identical individuals in the epidemic process with asymptotic 
probability $(l-1)/Nlp_l$,
since the element~$v$ has only $(l-1)$ half edges available to be matched, out of a total
number of half--edges from type~$l$ individuals that is still asymptotically $Nlp_l$.
Thus 
the mean number of such pairs that correspond to actual matches is asymptotically given by
\[
   \frac{(l-1)}{lp_l}\,\tbW\uid e^{\l u} \hz_l\, \int_0^\infty \l v e^{-\l v}\skK  \z_k(k-1)\{lp_l/m\}\, 
          (1-\Phi_k(v))\, G_{kl}(dv),
\]
and hence the probability that there is no such pair of any type~$l$, $1\le l\le d$, is
asymptotically given by $\exp\{-\tbW\uid e^{\l u} m_*\ut\}$, where
\eq\label{m-star-2-def-C}
   m_*\ut \Def \frac1{m}\skK\slK  \z_k (k-1)(l-1)\hz_l \int_0^\infty \l v e^{-\l v} (1-\Phi_k(v)) G_{kl}(dv)\,.
\en

These assertions, and the analogous assertions about the probability of two randomly
chosen individuals being infected by the initial individual, can be proved by the
methods introduced in Section~\ref{basic}, and lead to the following theorem.  
Here, $\ffptn$ denotes the
$\s$-algebra associated with the (forward) infection process until~$\nit$ infections have
occurred, and~$S_{Nl}(t)$ denotes the number of type~$l$ susceptibles at time~$t$.

\begin{theorem}\label{main-multitype-C}
Suppose that the forward branching process is supercritical.
Then there exists an event~$\tA_N \in \ffptn$ such that $\pr[\tA_N^c] \to 0$ as $N\to\infty$,
for which
\[
   \pr\Bigl[\sup_u |(Np_l)^{-1}S_{Nl}(\t_N + \l^{-1}\{\half\log N + u\}) -\hs_{l}(u)| > \e 
                   \Giv \ffptn\cap \tA_N \cap \{\t_N < \infty\} \Bigr] \ \to\ 0
\] 
as $N\to\infty$, where~$\hs_{l}$ is the decreasing function given by
\[
   \hs_{l}(u) \Def \tbaf\ul( e^u m_*\ut)  ,
\]
where $\tbaf\ul$  and~$m_*\ut$ are as defined above.  In particular, the total proportion of
susceptibles $N^{-1}\slK S_{Nl}(\t_N + \l^{-1}\{\half\log N + u\})$ is well approximated
by $\slK p_l \hs_l(u)$, uniformly in~$u$.
\end{theorem}

The general formulation above simplifies, if the distributions~$\Phi_k$ of infectious period
and~$G_{kl}$ of contact times are the same for all choices of the indices.  In this case,
the matrix~$\m(s)$ is given by
\[
    \m_{lk}(s) \Def (l-1)\{kp_k/m\}\,\int_0^\infty e^{-su}\,(1 - \Phi(u))G(du)
     \ =:\ U(s)(l-1)\{kp_k/m\},
\]
and is of rank one.  The positive eigenvalue is $U(s)m_{(2)}/m$, where
$m_{(2)} := \skK k(k-1)p_k$, the process is
supercritical if $m_{(2)}/m > 1/U(0)$, where $U(0) = \int_{(0,\infty)}(1 - \Phi(u))G(du)$,
and $\l$ is such that $U(\l) = m/m_{(2)}$.  The eigenvectors for the forward and
backward processes are equal, with $\z_i = \hz_i = ip_i/m$ and $\h_i = \heta_i = (i-1)m/m_{(2)}$. 
The quantities $m_*\ui$ and~$m_*\ut$ become
\[
   m_*\ui \Eq  m_0m_{(2)}/m \qquad\mbox{and}\qquad m_*\ut \Eq m_0 m_{(2)}^2/m^3,
\]
where $m_0 := \int_0^\infty \l v e^{-\l v} (1-\Phi(v)) G(dv)$.
The equations~\Ref{LT-vec-backward-C-back} can be much more neatly expressed, because
the functions~$\baf_0\ul$ are now the same for all~$l$, reflecting that the backward
process of half--edges is equivalent to a single--type branching process.
They reduce to the single equation
\eqa
  \baf_0(s) &=&  m^{-1}\skK kp_k 
    \Blb \int_{(0,\infty)} \{\baf_0(se^{-\l v})\}^{k-1}\,(1-\Phi(v))G(dv)  + 
      \Bl 1 - U(0) \Br \Brb \non\\
   &=& 1 - \int_{(0,\infty)}\{1 - m^{-1}g'(\baf_0(se^{-\l v}))\}(1-\Phi(v))G(dv),
     \label{LT-vec-backward-C-back-simple}
\ena
where $g(s) := \skK p_k s^k$, and the initial condition is
$(-D\baf_0)(0) = \heta_l/\{(l-1)m_*\ui\} = m/\{m_{(2)}m_*\ui\}$.
To express~$\hs_{l}(u) = \{\baf_0(e^{\l u} m_*\ut)\}^l$ more concisely, 
we write $h_s(u) := \baf_0(se^{\l u})$; then~\Ref{LT-vec-backward-C-back-simple} implies
that $h := h_s$ satisfies the equation
\eq\label{LT-vec-C-h}
   h(u) \Eq 1 - \int_{(0,\infty)}\{1 - m^{-1}g'(h(u-v))\}(1-\Phi(v))G(dv),
\en
with initial condition $\lim_{u\to-\infty}\{e^{-\l u}(Dh)(u)\} = \l s$;
so $\hs_1(u) = h_{m_*\ut}(u)$ satisfies~\Ref{LT-vec-C-h} with 
$\lim_{u\to-\infty}\{e^{-\l u}(D\hs_1)(u)\} = \l m_*\ut$, and $\hs_l = (\hs_1)^l$.

In the case of the Volz~(2008) model, there is further simplification, because of
the explicit forms $\Phi(v) = 1 - e^{-\b v}$ and $G(dv) = \a e^{-\a v} dv$.  In this case,
a deterministic law of large numbers starting with an asymptotically positive initial proportion
of infectious individuals was established by Decreusefond {\it et al.\/} (2012).  With
such an initial condition, the randomness inherent in the initial stages of development,
reflected in the presence of~$\t_N$ in the statement of Theorem~\ref{main-multitype-C},
plays no significant part.  In the Volz setting,
explicit formulae for $\l = \a m_{(2)}/m - \b$ and $m_0 = \l\a(\l+\a+\b)^{-2}$ can be written
down, and equation~\Ref{LT-vec-C-h} can be expressed as
\eqs
   h(u) &=& 1 - \int_{(0,\infty)}\{1 - m^{-1}g'(h(u-v))\}\a e^{-(\a+\b)v}\,dv \non \\
        &=& \frac\b{\a+\b} + \frac1m \int_{-\infty}^u g'(h(w)) \a e^{-(\a+\b)(u-w)}\,dw\,.
\ens
Differentiation with respect to~$u$ then yields the following autonomous differential 
equation for~$h = h(t)$: 
\eq\label{LT-vec-C-h-new}
  \frac{dh}{dt} \Eq \frac\a{m} g'(h) - (\a+\b)h + \b \Eq (\a+\b)(\tf(h) - h),
\en
where~$\tf(s)$ is the probability generating function $(\a g'(s) + \b)/(\a+\b)$.
In particular, it follows that $h(\infty)$ is the solution~$\tq$ smaller than~$1$ to the
equation $\tf(s) = s$, and hence that the asymptotic final proportion of susceptible 
individuals at the end of a large outbreak is given by $g(\tq)$.

\medskip
\begin{remark}
Volz~(2008) expresses the equations for the development of the epidemic
as the solutions to a system of three coupled differential equations for the variables
$h$, $p_I$ and~$p_S$:
\eqs
   \frac{dh}{dt} &=& - \a p_I h;
   \qquad\frac{dp_S}{dt} \Eq \a p_Sp_I \Bl 1 - \frac{hg''(h)}{g'(h)} \Br; \\
   \frac{dp_I}{dt} &=& \a p_Sp_I \frac{hg''(h)}{g'(h)}
             - \a p_I(1-p_I) - \b p_I.
\ens
It is not difficult to see that their solution is given in terms of the solution~$h$ 
to~\Ref{LT-vec-C-h-new} by $p_I = 1 - \frac{g'(h)}{mh}  + \frac{\b}\a\{1 - \frac1h\}$ and 
$p_S = \frac{g'(h)}{mh}$.  The first equation is clearly satisfied, and the second follows
by differentiating the formula for~$p_S$ and using the first equation to re-express $\frac{dh}{dt}$.
Then the sum of the second and third equations is satisfied by differentiating the
expression for $p_S + p_I$, and then using it once more to re-express $(1 - 1/h)$.
\end{remark}

\medskip
\begin{remark}
   Although the asymptotics carried out in this section are not applicable to
that case, a reasonably general Kermack--McKendrick epidemic also fits into this
epidemic model, by taking $p_{N-1} = 1$ and by replacing $G(du)$ by $(N-1)^{-1}G(du)$;
in the notation of Section~\ref{intro}, we would have~$\b(u)du$ replacing $(1-\Phi(u))G(du)$.
This leads formally to equations determining the development of the epidemic
which are asymptotically equivalent, for large~$N$, to those established in 
Section~\ref{basic}.  For instance, \Ref{LT-vec-C-h} becomes
\[
     h(u) \Eq 1 - \frac1{N-1}\int_0^\infty \{1 - h^{N-2}(u-v)\}\b(u)\,du,
\]
so that, writing $\hs(u)$ for $h(u)^{N-1}$, we obtain
\[
     \hs(u) \sim \exp\Blb - \int_0^\infty \{1 - \hs(u-v)\}\b(u)\,du \Brb,
\]
which is just~\Ref{deterministic}.
\end{remark}

\section*{Acknowledgement}
We are grateful to Peter Jagers, Hans Heesterbeek and Odo Diekmann for a number of
helpful discussions.
ADB thanks  the mathematics departments of the University of Melbourne, Monash University and the
University of Queensland, and the Institute for Mathematical Sciences at the National University
of Singapore, for
their kind hospitality while part of the work was undertaken. 

\end{document}